\documentclass[12pt, reqno]{amsart}

\usepackage{etex} 
\usepackage{amsmath, amsthm, amssymb, mathtools}  
\usepackage{graphicx, color, bm, url}  
\usepackage{caption, subcaption}  
\usepackage{faktor}  
\usepackage{enumerate, enumitem}  
\usepackage{verbatim}  
\usepackage{chemfig, chemnum}  
\usepackage{tikz, tikz-cd}  
\usepackage{diagbox}  
\usepackage{float}  
\usepackage{array}  

\usepackage[english]{babel}
\usepackage[utf8]{inputenc}  
\usepackage[T1]{fontenc}  

\usepackage[colorlinks=true, allcolors=blue]{hyperref}
\usepackage{cite}  
\usepackage[capitalise]{cleveref}  

\usepackage[top=30mm, bottom=30mm, left=25mm, right=25mm]{geometry}
\usepackage{ytableau}

\setlist[1]{itemsep=3pt}  

\usepackage{setspace}

\newtheorem{theorem}{Theorem}[section]
\newtheorem{lemma}[theorem]{Lemma}
\newtheorem{proposition}[theorem]{Proposition}
\newtheorem{corollary}[theorem]{Corollary}

\theoremstyle{definition}
\newtheorem{definition}[theorem]{Definition}
\newtheorem{example}[theorem]{Example}

\newtheorem{conjecture}[theorem]{Conjecture}
\theoremstyle{remark}
\newtheorem{remark}[theorem]{Remark}

\numberwithin{equation}{section}

\newcommand{\beq}{\begin{equation}}
\newcommand{\eeq}{\end{equation}}

\DeclareMathOperator{\KP}{K\mathrm{P}}

\DeclareMathOperator{\B}{B}
\DeclareMathOperator{\id}{id}
\DeclareMathOperator{\val}{val}
\DeclareMathOperator{\Aut}{Aut}
\DeclareMathOperator{\GL}{GL}

\DeclareMathOperator{\End}{\mathrm{End}}
\DeclareMathOperator{\Hom}{Hom}
\DeclareMathOperator{\diag}{\mathrm{diag}}

\DeclareMathOperator{\Fix}{Fix}
\DeclareMathOperator{\carac}{char}
\DeclareMathOperator{\SO}{SO}

\DeclareMathOperator{\Span}{\mathrm{span}}
\DeclareMathOperator{\Frac}{Frac}
\DeclareMathOperator{\image}{im}
\DeclareMathOperator{\supp}{supp}
\DeclareMathOperator{\Stab}{\mathrm{Stab}}
\DeclareMathOperator{\PGL}{PGL}

\renewcommand{\P}{\mathbb{P}}
\newcommand{\N}{\mathbb{N}}
\newcommand{\R}{\mathbb{R}}
\newcommand{\Q}{\mathbb{Q}}
\newcommand{\F}{\mathbb{F}}
\newcommand{\C}{\mathbb{C} }
\newcommand{\Z}{\mathbb{Z}}
\newcommand{\Ocal}{\mathcal{O}}
\newcommand{\Hcal}{\mathcal{H}}

\newcommand{\Bcal}{\mathcal{B}}

\begin{document}

\title[]{$\GL(n,\Z_p)$-invariant Gaussian measures on the space of $p$-adic polynomials}
\author[Y.~El~Maazouz]{Yassine~El~Maazouz}
\address{California Institute of Technology, Pasadena, California.}
\email{maazouz@caltech.edu}

\author[A.~Lerario]{Antonio~Lerario}
\address{SISSA, Trieste, Italy.}
\email{lerario@sissa.it}




\dedicatory{}


    \begin{abstract}
        We prove that if $p>d$ there is a unique gaussian distribution (in the sense of Evans \cite{evansLocalFields2001}) on the space $\Q_p[x_1, \ldots, x_n]_{(d)}$ which is invariant under the action of $\GL(n, \Z_p)$ by change of variables. 
        
        This gives the nonarchimedean counterpart of Kostlan's Theorem \cite{Kostlan93} on the classification of orthogonally (respectively unitarily) invariant gaussian measures on the space $\R[x_1, \ldots, x_n]_{(d)}$ (respectively $\C[x_1, \ldots, x_n]_{(d)}$).
             
         More generally, if  $V$ is an $n$--dimensional vector space over a nonarchimedean local field $K$ with ring of integers $R$, and  if $\lambda$ is a partition of an integer $d$, we study the problem of determining the invariant lattices in the Schur module $S_\lambda(V)$ under the action of the group $\GL(n,R)$.
    \end{abstract}
\maketitle

\section{Introduction}

\subsection{Probabilistic Algebraic Geometry.}\label{sec:probabilityintro}
	 
	  In the last years there has been an increasing interest into the statistical behaviour of algebraic sets over non-algebraically closed fields. 
	 When the notion of ``generic'' is no longer available, one seeks for a ``random'' study of the objects of interest. 
	 For instance, if $K$ is a field, once a probability distribution is put on the space $K[x_1, \ldots, x_n]_{(d)}$ of homogeneous polynomials, one can study expected properties of their zero sets in projective space $K\mathrm{P}^{n-1}$ (e.g. the expected number of solutions of systems of random equations), see for instance \cite{EdelmanKostlan95,EKS, E21, ET19, EK22, Letwo, GaWe1, GaWe2, GaWe3, ShSm1, ShSm3, Ko2000, LeLu:gap, NazarovSodin1, Sarnak, shsm} for the real and the complex case and \cite{BCFG22,Evans, AL20,Caruso22,KulkarniLerario, shmueli} for the $p$--adic case.

	 The appropriate choice of a norm on $K^n$ (using a scalar product when $K=\R$, a hermitian product when $K=\C$ or a non-archimedean norm when $K=\Q_p$) induces a metric structure on $K\mathrm{P}^{n-1}$ and the group $\mathrm{Iso}(K^n)\subset \GL(K^n)$ of linear, norm preserving transformations, acts by isometries on the projective space.
     In this context it is natural to put a probability distribution on the space of polynomials $K[x_1, \ldots, x_n]_{(d)}$ which is invariant under the action of $\mathrm{Iso}(K^n)$ induced by change of variables -- there should be no preferred points or directions in projective space.
     An interesting problem is  therefore to find (and possibly classify) all the probability distributions on $K[x_1, \ldots, x_n]_{(d)}$ having this property.
	 
	 When $K=\R, \C$, and in the case of nondegenerate \emph{gaussian} distributions, this problem is equivalent to the problem of finding scalar products (hermitian structures in the complex case) on $K[x_1, \ldots, x_n]_{(d)}$ which are preserved by orthogonal change of variables, and it was solved by Kostlan \cite{Kostlan93} using representation theory.
	 
	  More precisely, if $K=\C$, then $\mathrm{Iso}(\C^n)\simeq U(n,\C)$ and, since the representation
	 $$\rho_{n,d}^\C:U(n, \C)\to \GL(\C[x_1, \ldots, x_n]_{(d)})$$ is irreducible, by Schur's Lemma there is a unique such hermitian structure (up to multiples). This is called the Bombieri--Weyl hermitian structure.
	 
	 If $K=\R$, then $\mathrm{Iso}(\R^n)\simeq O(n, \R)$ and the representation 
	 \[\rho_{n,d}^\R:O(n, \R)\to \GL(\R[x_1, \ldots, x_n]_{(d)})\] \emph{is not} irreducible. Its irreducible summands are spaces of spherical harmonics and there is a whole $\lfloor \frac{d}{2}\rfloor$-dimensional family of scalar products having the desired property, see \cite{Kostlan93}.

	 	 If  $K=\mathbb{Q}_p$, then $\mathrm{Iso}(\Q_p^n)\simeq\GL(n, \Z_p)$ (see \cite[Theorem 2.4]{EvansRaban19}). Moreover, Evans \cite{Evans88,Evan88_2,Evans90,Evans93,evansLocalFields2001,Evans07}  introduced a notion of  \emph{gaussian} distribution on a $p$-adic vector space, which essentially corresponds to a choice of lattice in the vector space (in the same way as real gaussian structures corresponds to scalar products, see Section \ref{sec:probability}). Using this correspondence, in this context the problem of determining the invariant distributions can be formulated as: find all lattices in $\Q_p[x_1, \ldots, x_n]_{(d)}$ which are invariant under the action of $\GL(n, \Z_p)$ by change of variables. As in the real case, these invariant distributions have nice properties in terms of computing expectation of number of zeroes of system of random equations, see Section \ref{sec:expected}.

		 Theorem \ref{thm:FixIsOnePointintro} implies that, if  $p>d$, there is only one such lattice (up to scaling) and therefore only one probability distribution on $\Q_p[x_1, \ldots, x_n]_{(d)}$ with the required property, up to scaling (see Section \ref{sec:probability} for more details).
		 
		 \begin{theorem}\label{theorem:probability}
		There are finitely many gaussian distributions  on $\Q_p[x_1, \ldots, x_n]_{(d)}$ (up to multiples) which are invariant under the action of 
		 $\GL(n, \Z_p)$ by change of variables. Moreover, if $p>d$, there is only one such distribution (up to multiples), which is given by the uniform measure\footnote{This means that a random polynomial sampled from this distribution can be written as 
		 $$f(x)=\sum_{|\alpha|=d}\zeta_{\alpha}x_1^{\alpha_1}\cdots x_n^{\alpha_n},$$ where the $\zeta_\alpha$ are independent uniform random variables in $\Z_p.$} on $\Z_p[x_1, \ldots, x_n]_{(d)}$.		 
		 \end{theorem}
          The irreducibility of $\rho_{n,d}:\GL(n, \Z_p)\to \GL(\Q_p[x_1, \ldots, x_n]_{(d)})$ (Theorem \ref{thm:repIsIrred} below) plays a role for the uniqueness, but in a different way (compared to the complex case). More generally, we have the following result, proved in Section \ref{sec:proofofcoro}, which could be of some interest, for instance, for the study of random tensors over nonarchimedean local fields.
    
            \begin{example}[$n=2, d=2$] The group $\GL(2,\Z_p)$ acts on the vector space of degree $2$ polynomials $\Q_p[x,y]_{(2)} = \Span_{\Q_p}( x^2 , xy, y^2)$ as follows:
            \[
                (g \cdot P)(x,y) := P(ax + cy, bx + dy), \quad \text{for } g \in \GL(2, \Z_p) \text{ and } P \in \Q_p[x,y]_{(2)}.
            \]
            When $p > 2$ is an odd prime number, the only Gaussian probability measures on $\Q_p[x,y]_{(2)}$ that are invariant under the action of $\GL(2,\Z_p)$ are the measures of the following form
            \[
                \bm{\zeta} =  p^{r} \Big( \zeta_1 \ x^2 + \zeta_2 \  xy + \zeta_3  \ y^2 \Big),
            \]
            where $\zeta_1, \zeta_2, \zeta_3$ are i.i.d uniformly distributed random variables in $\Z_p$ and $r \in \Z \cup \{+\infty\}$. The condition $p > 2$ in Theorem \ref{theorem:probability} is crucial, see Example \ref{exa:counterExample}. 
            \end{example}
	 
    	\begin{proposition}\label{prop:InvariantMeasures}
    		    Suppose $K$ is a non-archimedean local field and assume that $\lambda$ is $\carac(K)$-core.
    		    Then there are only finitely many Gaussian distributions on $S_\lambda(V)$ (up to scaling) that are invariant under the action of $\GL(n,R)$ through $\rho_{n,\lambda}$. 
    		    Moreover, if $\lambda$ is also $p $-core, where $p$ is the residue characteristic of $K$, there is only one such measure (up to scaling).
     	\end{proposition}

       Recall that the hook length of a box in the Young diagram of a partition $\lambda$ is obtained by adding $1$ to the number of boxes below and to the right of the box in question. If a prime $\ell$ does not divide any of the hook lengths of $\lambda$, the partition is said to be \emph{$\ell$-core} (by convention every partition is $0$-core). In particular the trivial partition $\lambda = (d)$ is $p$-core if an only if $p>d$.

\subsection{Invariant lattices.} Let us discuss our problem using the notion of lattices and in a general framework.

	 Let a $K$ be a discretely valued non-archimedean field and $R$ its valuation ring. Given a vector space $V$ of dimension $n$ over $K$ (identified with $K^n$ through the choice of a basis) and a partition $\lambda$ of an integer $d$, we are interested in determining lattices in the Schur module $S_\lambda(V)$ which are invariant under the action of $\GL(n, R)$ given by the representation
	 \begin{equation}\rho_{n, \lambda}:\GL(n, R)\to \GL(S_{\lambda}(V)).\end{equation}      
    
    In the case  $K=\mathbb{Q}_p$, $R=\mathbb{Z}_p$ and $\lambda=(d)$, we have $S_{\lambda}(V) \simeq \mathbb{Q}_p[x_1, \ldots, x_n]_{(d)}$, the space of homogenous polynomials in $n$ variables of degree $d$, with the action of $\GL(n, \mathbb{Z}_p)$ being the action by change of variables. As we already mentioned, the problem of determining lattices in $\mathbb{Q}_p[x_1, \ldots, x_n]_{(d)}$ which are invariant under this action is equivalent to the problem of determining gaussian probability distributions on $\mathbb{Q}_p[x_1, \ldots, x_n]_{(d)}$.
    
    With this notation, the next theorem is our main technical result. 
    
    \begin{theorem}\label{thm:FixIsOnePointintro}
         Denote by $\ell$ the characteristic of $K$ and by $p$ its residue characteristic. If the partition $\lambda$ is $\ell$-core, then there are only finitely many $\GL(n, R)$-invariant lattices in $S_\lambda(V)$ up to scaling. If moreover $\lambda$ is $p$-core, then:
        \[
                \Fix(\rho_{n,\lambda}^{\Bcal}) = \left \{ [S_\lambda(L_0)] \right \},
        \]
        with $L_0 \coloneqq R \cdot x_1 \oplus \dots \oplus R \cdot x_n$ is the unique $\GL(n, R)$-invariant lattice in $V$, up to scaling.
    \end{theorem}


\subsection{Group actions on the Bruhat-Tits Buildings.} Our main result can be restated more precisely (see Theorem \ref{thm:FixIsOnePoint}) using the language of \emph{Buildings}. These are combinatorial and geometric structures that generalize certain aspects of Riemannian symmetric spaces, see \cite{Ji09,Rousseau09}.
    
    In our context, observe that lattices in $S_{\lambda}(V)$ (up to homothety) are points in the \emph{Bruhat-Tits building} $\Bcal_{n,\lambda}$ for $\PGL(S_\lambda(V))$. This is an infinite simplicial complex whose $0$-simplices are the homothety classes of lattices in $S_\lambda(V)$, see Section \ref{sec:buildings}. The action of $\GL(n, R)$ on $S_{\lambda}(V)$ induces an action on the building $\Bcal_{n,\lambda}$ i.e. we have a group homomorphism
			\[
				\rho_{n,\lambda}^{\Bcal} \colon \GL(n, R) \to \Aut(\Bcal_{n,\lambda}),
			\]
	sending $\GL(n, R)$ to the group of automorphisms of the building. With this notation, Theorem \ref{thm:FixIsOnePoint} below states that if $\lambda$ is $\carac(K)$-core, the set of fixed points of this action is a non-empty finite convex set in the building $\Bcal_{n,\lambda}$, reduced to one point if $\lambda$ is $p$-core, where $p$ is the residue characteristic of  $K$.

	\subsection{Structure of the paper}
         This article is organized as follows. In Section \ref{sec:probability}, we explain the probabilistic motivation behind this work. We collect some background and preliminary results in Section \ref{sec:Background}, and prove our main results in Section \ref{sec:proofs}. In Section \ref{sec:slices} we give a way to compute  apartment slices of $\Fix(\rho_{n,\lambda})$. In Section \ref{sec:expected} we discuss $p$-adic analogues of the work of Edelman, Kostlan, Shub and Smale on system of random equations (this section serves as a motivation for the results). Finally, we discuss some open questions and final remarks in Section \ref{sec:concluding_remarks}.
        
        \subsection*{Acknowledgements.}
        This work started during Y.E's visit to the MPI MiS and  SISSA.
        We thank the two aforementioned institutions for their generous hospitality. The authors thank Kiran Kedlaya, Vera Serganova, Bernd Sturmfels and Annette Werner for helpful discussions.
    
    \section{Probabilistic motivation} \label{sec:probability}

            We recall in this section some standard notions from probability, to help the reader to put our results into context. The reader may skip this section and come back to it when necessary.
            
            \subsection{Gaussian distributions on real vector spaces.} 
            We start by recalling the notion of gaussian distributions on a real vector space. First, the \emph{standard gaussian distribution} on $\R^m$ is the probability measure $\gamma_m$ on $\R^m$ defined for every Borel set $U\subseteq \R^m$ by
            \[\gamma_m(U):=\frac{1}{(2\pi)^{\frac{m}{2}}}\int_U e^{-\frac{\|x\|^2}{2}}dx,\]
            where $\|x\| \coloneqq \sqrt{|x_1|^2 + \dots + |x_n|^2}$ denotes the Euclidean norm of $x=(x_1, \ldots, x_m)$.
            In practice, if $\{e_1, \ldots, e_m\}$ is an orthonormal basis for the Euclidean norm and $\xi_1, \ldots, \xi_m$ are independent, standard gaussians (i.e. $\P(\xi_j\leq t)=\gamma_1(-\infty, t)$ for every $j=1, \ldots, m$), sampling from the standard gaussian distribution on $\R^m$ is equivalent to picking an element $\xi \in \R^m$  at random by writing it as a random linear combination
            \[\xi=\xi_1 e_1+\cdots +\xi_m e_m.\]
            
            More generally, a nondegenerate, centered, gaussian distribution on a real finite dimensional vector space $V$ is given by assigning a surjective linear map $T:\R^m\to V$ and defining $\P(W):=\gamma_m(T^{-1}(W))$ for every Borel set $W\subseteq V$. (This measure is denoted by $T_\#\gamma_m$ and called the \emph{push-forward} measure.) Writing $\R^m=\ker(T)\oplus \ker(T)^\perp$, one gets a linear isomorphism $V\simeq \ker(T)^\perp$ and an induced scalar product on $V$. As above, if $\{v_1, \ldots, v_n\}$ denotes an orthonormal basis for $V$, this construction is equivalent to define a random element in $V$ by $\xi_1v_1+\cdots+\xi_n v_n$, with the $\xi_j$'s standard, independent gaussians.
            Viceversa, given a scalar product on $V$ and an orthonormal basis $\{v_1, \ldots, v_n\}$ for it, putting $\xi:=\xi_1v_1+\cdots+\xi_n v_n$, where the $\xi_i$ are standard, independent gaussians, defines a random variables with values in $V$ and the induced probability distribution is nondegenerate, centered and gaussian (as a map $T:\R^n\to V$ in this case we could simply take $T(e_i):=v_i$).
            
            From this we see that the theory of nondegenerate, centered, gaussian distributions on $V$ is equivalent to the theory of nondegenerate positive definite quadratic forms on it (i.e. scalar products). In particular, once a representation $\rho:G\to \GL(V)$ is given, asking for the nondegenerate, centered gaussian probability distributions $\gamma$ on $V$ which are $\rho$-invariant (i.e. such that $\rho(g)_{\#}\gamma=\gamma$ for all $g\in G$) is equivalent to ask for the scalar products $\langle\cdot,\cdot\rangle$ on $V$ which are $\rho$-invariant (i.e. such that $\langle v_1, v_2\rangle=\langle \rho(g)v_1, \rho(g)v_2\rangle$ for all $v_1, v_2\in V$ and $g\in G$). 
            
            In particular, if the representation $\rho$ is irreducible, Schur's Lemma implies that there is only one such invariant scalar product (up to multiples) and, consequently, only one invariant gaussian distribution (up to scaling). If the representation \emph{is not} irreducible, this is no longer true as it happens for the case of real polynomials.
            
            \subsection{Invariant gaussian distributions on the space of real and complex polynomials.} 
            Let now $\rho_{n,d}^\C :U(n, \C) \to \GL(\C[x_1, \ldots, x_n]_{(d)})$ be the representation given by change of variables:
           \[\rho_{n,d}^\C(g)(P):=P\circ g^{-1}.\]
           This representation is complex irreducible and, consequently, there is only one $\rho_{n,d}^\C$-invariant hermitian structure (up to multiples) on $\C[x_1, \ldots, x_n]_{(d)}$. It is called the \emph{Bombieri-Weyl} hermitian structure. A hermitian orthonormal basis for it is given by the monomials
	 \begin{equation}\label{eq:BW}\left\{\left(\frac{d!}{\alpha_1!\cdots \alpha_n!}\right)^{1/2}x_1^{\alpha_1}\cdots x_n^{\alpha_n}\right\}_{\alpha_1+\cdots+\alpha_n=d}.
	 \end{equation}
	 The Bombieri-Weyl scalar product and the corresponding gaussian distribution have been widely used, see  \cite{EdelmanKostlan95, shsm, EKS, Ko2000, ShSm3, ShSm1, GaWe1, GaWe3, GaWe2, NazarovSodin1, Sarnak, Letwo, LeLu:gap}.
	 
	      Notice that the basis in \eqref{eq:BW} is real and, since $O(n, \R)\subset U(n, \R)$, the Bombieri--Weyl scalar product restricts to a scalar product on $\R[x_1,\ldots, x_n]_{(d)}$ which is invariant under $O(n, \R)$. In other words, denoting by $\rho_{n,d}^\R : O(n, \R) \to \GL(\R[x_1,\ldots, x_n]_{(d)})$ the representation by change of variables,  the Bombieri-Weyl scalar product is $\rho_{n,d}^\R$-invariant. However, since $\rho_{n,d}^\R$ \emph{is not} real irreducible, there are invariant scalar products which are not multiples of this one. The classification of such scalar products has been done by Kostlan in \cite{Kostlan93} using the theory of spherical harmonics, as follows. Denoting by $\Hcal_{n,\ell}\subset \R[x_1, \ldots, x_n]_{(\ell)}$ the space of harmonic polynomials, we have a decomposition
	      $$\R[x_1, \ldots, x_n]_{(d)}\simeq \bigoplus_{d-\ell \text{ even}}\|x\|^{d-\ell}\cdot \Hcal_{n, \ell}.$$
          The spaces $\|x\|^{d-\ell}\cdot \Hcal_{n, \ell}$ are precisely the irreducible summands of $\rho_{n,d}^\R$ and are isomorphic to the spaces of spherical harmonics.   Schur's Lemma implies that there is only one invariant scalar product on each of them, up to multiples, and an invariant scalar product on $\R[x_1, \ldots, x_n]_{(d)}$ is obtained by scaling these scalar products separately. For instance, the $L^2(S^{n-1})$ scalar product is $\rho_{n,d}^\R$-invariant but is not a multiple of the Bombieri-Weyl one (see \cite{LerarioFLL}).      
    
          \subsection{Gaussian distributions on \texorpdfstring{$p$}{}-adic spaces.}
          Evans \cite{Evans88,Evan88_2,Evans90,Evans93,evansLocalFields2001,Evans07} introduced the notion of gaussian distribution on a $p$-adic vector space, as follows. To start with, one denotes by $\zeta_1$ the uniform probability measure on the compact topological group $\Z_p$ and by $\zeta_m:=\zeta_1\times \cdots\times  \zeta_1$ the product measure on $\Z_p^m$ (these are just the normalized Haar measures). The measure $\zeta_m$ is called the \emph{standard $p$-adic gaussian measure}. Then, if $V\simeq \Q_p^n$ is a $p$-adic vector space, a $p$-adic gaussian measure is defined by assigning a surjective linear map $T:\Q_p^m\to V$ and considering the pushforward measure $T_{\#}\zeta_m$ (notice the analogy with the real construction).
          
          Here the theory of nondegenerate $p$-adic gaussian measures on $V\simeq \Q_p^n$ is equivalent to the theory of full dimensional lattices in $V$: the image of $\Z_p^m$ under $T$ is a lattice $L:=T(\Z_p^m)$ in $V$ and every lattice arise in this way. This lattice is the support of the measure $T_{\#}\zeta_m$. Given a lattice $L\subset V\simeq \Q_p^n$, one defines a $p$-adic gaussian distribution fixing a basis $\{v_1, \ldots, v_n\}$ for $L$ and setting $T(e_j)=v_j$, as above (here $\{e_1, \ldots, e_n\}$ is the standard basis). A random element from this distribution is obtained by setting $\xi:=\xi_1 v_1+\cdots +\xi_n v_n$, where the $\xi_j$'s are independent uniform random variables on $\Z_p$. Gaussian distributions on a discretely valued non-archimedean field $K$ are constructed similarly.
          
          Given a representation $\rho:G\to \GL(V)$, it is then natural to ask for lattices which are invariant under this representation (or equivalently invariant under the action of the ring and $R$-module $H_\rho \coloneqq \Span_R(\image(\rho)) \subset \End_K(V)$); they correspond to $\rho$-invariant gaussian distributions. Unlike the real and complex setting, in the non-archimedean case, even if $\rho$ is irreducible, Schur's Lemma cannot be used to conclude uniqueness in this context. That is because non-degenerate gaussian measures are in a one-to-one correspondence with lattices (instead of positive non-degenerate quadratic forms as in the real case).

         \subsection{Invariant gaussian distributions on the space of \texorpdfstring{$p$-}{}adic polynomials.}         
          Evans \cite{Evans} defines a probability distribution on the space of polynomials considering the random polynomial
    	 \begin{equation}
    	    \zeta(y):= \zeta_0 + \zeta_1 \binom{y}{1} + \cdots +\zeta_d \binom{y}{d},
    	 \end{equation}
    	 where $\zeta_0, \ldots, \zeta_d$ are independent and uniformly distributed in $\Z_p$. Since the aforementioned seminal work of Evans a couple of decades ago, probabilistic problems over non-archimedean local fields have been gaining interest in the recent years \cite{ET19,E21,AL20,BCFG22,Caruso22,EK22}. In many of these problems, it is important that the probability measure is invariant under certain symmetries. Homogenizing the above polynomial, one gets a probability distribution on $\Q_p[x, y]_{(d)}$ which \emph{is not} invariant under the action of $\GL(2, \Z_p)$ by change of variables. Here $\GL(n, \Z_p)$ can be seen as the group of isometries of projective space $\Q\mathrm{P}^{n-1}$ (see \cite{KulkarniLerario}), and it is natural to ask for a probability distribution on $\Q_p[x_1, \ldots, x_n]_{(d)}$ for which there are no preferred points or direction for the zero sets of polynomials in projective space, as we did for the real case.

    	 In \cite{KulkarniLerario} the authors proposed an alternative model, defining a random polynomial $\zeta\in \Q_p[x_1, \ldots, x_n]_{(d)}$ as
    	 \begin{equation}\label{eq:inv}\zeta(x):=\sum_{|\alpha|=\alpha_1+\cdots +\alpha_n=d}\zeta_{\alpha}x_1^{\alpha_1}\cdots x_n^{\alpha_n},
    	 \end{equation}
    	 where $\{\xi_\alpha\}_{|\alpha|=d}$ is a family of independent random variables, each of them uniformly distributed in $\Z_p$. This gaussian distribution corresponds to the lattice $L\subset \Q_p[x_1, \ldots, x_n]_{(d)}$ spanned by the standard monomial basis. This has been used, for example, in \cite{AL20} to compute the expected number of $k$--planes on a random complete intersection (e.g. the number of lines on a random cubic surface).
    	 	
    	 The probability distribution induced by \eqref{eq:inv} on $\Q_p[x_1, \ldots, x_n]_{(d)}$ is invariant under the action of $\GL(n, \Z_p)$ by change of variables (see \cite{KulkarniLerario}). It is natural to ask whether there are other Gaussian distributions on $\Q_p[x_1, \ldots, x_n]_{(d)}$ which have this property. Theorem \ref{thm:FixIsOnePointintro} implies that, when $p>d$, \eqref{eq:inv} is the only distribution with this property (up to scaling). In fact, this corresponds to the case $\lambda = (d)$ with the Young diagram
        \begin{center}
            $\lambda = 
            \ytableausetup{centertableaux}
            \begin{ytableau}
                    1 & 2 & \cdots  & d \\
            \end{ytableau}$
        \end{center}
        More generally, the same result is true for a  discretely valued non-archimedean field $K$: Proposition \ref{prop:InvariantMeasures} states that if $\lambda$ is $p$-core, with $p$ the residue characteristic, then there is a unique Gaussian measure on the space $S_\lambda(V) \cong K[x_1, \dots, x_n]_{(d)}$ of homogeneous polynomials of degree $d$ in $n$ variables that is invariant under linear change of variables in $\GL(n,R)$.

    \section{Background and preliminaries}\label{sec:Background}

	       \subsection{The group \texorpdfstring{$\GL(n, R)$}~.} 	To set things up, let $K$ be a non-archimedean discretely valued field with normalized valuation $\val \colon K^{\times} \twoheadrightarrow \Z$. We denote by $R$ its valuation ring and we fix a uniformizer $\varpi$ of $K$. We denote by $k$ the residue field $R / \varpi R$, and $\ell, p$ the respective characteristics \footnote{When $\ell > 0$ we necessarily have $p = \ell$.} of the fields $K$ and $k$. The non-archimedean valuation $\val$ induces an ultrametric absolute value $|\cdot|$ on $K$ as follows:
			$$ |x| \coloneqq \begin{cases}  p^{- \val(x)}, 	 \quad  \quad  \text{if } p > 0 \\ e^{-\val(x)}, \quad \quad \text{otherwise} \end{cases}. $$
			There is a natural notion of non-archimedean orthogonality (see \cite[Section 3]{evansLocalFields2001} for more details) for which the analogue of the group of orthogonal $n \times n$ matrices $O(n,\R)$ in the real setting is the group
			$$	\GL(n,R) \coloneqq \left\{ g \in \GL(n,K) \colon g, g^{-1} \in R^{n \times n} \right\}, $$
			see \cite[Theorem 2.4]{EvansRaban19}; this group is a totally disconnected compact topological group.
			
			\subsection{The Schur functor.} Let $V$ be an $n$-dimensional $K$-vector space\footnote{Note that same construction can be carried out more generally when $V$ is an $R$-module.} spanned by a basis $x_1, \dots, x_n$. There exists a natural linear right action of  the group $\GL(n, R)$ on $V$ as follows
			\[
					 x_i \cdot g = \sum_{j = 1}^{n} g_{ij} \ x_j,  \quad g = \left( g_{ij} \right)_{1 \leq i,j \leq n} \in \GL(n, R).
			\]
			This clearly defines a faithful representation $\rho_n \colon \GL(n, R) \to \GL(V)$. Let $d \geq 1$ and $\lambda \vdash d$ a partition of $d$. The \emph{Schur functor} $S_\lambda$ maps the representation $\rho_n$ to the representation $\rho_{n,\lambda}$ acting on the \emph{Weyl module} $S_{\lambda}(V)$ defined as
			\[
				S_{\lambda}(V) \coloneqq c_\lambda \cdot V^{\otimes d},
			\]
			where $c_\lambda$ is the \emph{Young symmetrizer}. For a detailed construction of the Weyl module $S_\lambda(V)$ we refer the reader to \cite[Section 6.1]{Fulton91}, \cite[Chapter 8]{Fulton97} or \cite[Chapter 6]{GrusonSerganova18}.

			Theorem 1 in section 8.1 of \cite{Fulton97} gives a basis of $S_\lambda(V)$ whose elements are indexed by the Young tableaux $T$ obtained by filling the Young diagram of $\lambda$ with entries in $\{1,\dots,d\}$. Choosing this basis as the \emph{standard basis} of $S_\lambda(V)$ we get a group homomorphism
			\[
			    \rho_{n,\lambda} \colon \GL(n,K) \to \GL(N,K),
			\]
            with $N \coloneqq \dim_K(S_\lambda(V))$.
			\begin{remark}
			    For geometric minded readers, the representation $\rho_{n,\lambda}$ is a homomorphism
			    $$ \rho_{n, \lambda} \colon \GL(n,\cdot) \to \GL(N, \cdot) $$
			    of group schemes over $\Z$.
			\end{remark}
			
			\subsection{Lattices.}			
			\emph{Lattices} in $S_{\lambda}(V)$ are full rank $R$-submodules of $S_{\lambda}(V)$. The representation $\rho_{n,\lambda}$ defines a natural action of $\GL(n, R)$ on lattices of $S_{\lambda}(V)$ as follows:
			\[
		 		g \cdot L  \coloneqq \{ g \cdot x  \colon x \in L \},  \quad \text{ as an $R$-module}.
			\]
			In this article, we are interested in determining the lattices that are invariant under this action i.e. lattices $L$ such that
			\[
				g \cdot L = L, \quad \text{for all } g \in \GL(n, R).
			\]
			(Obviously, if $L$ is $\rho_{n, \lambda}$-invariant, then $a L$ is $\rho_{n, \lambda}$-invariant for any $a \in K^\times$. So, to be more precise, we are interested in $\rho_{n, \lambda}$-invariant \emph{homothety classes} of lattices.)
			
			A central object for our study is the $R$-submodule of $\End_K(S_{\lambda}(V))$ generated by the image of $\rho_{n, \lambda}$, which we denote by $H_{n ,\lambda}$ i.e.
    	\[
    		\label{eq:defH}H_{n ,\lambda} = \Span_R(\image(\rho_{n,\lambda})).
    	\]
    	The $R$-module $H_{n, \lambda}$ is also a subring of $\End_K(S_{\lambda}(V))$. Given a lattice $L$ in $S_{\lambda}(V)$, the action of $\rho_{n,\lambda}$ defines a new lattice $H_{n,\lambda} \cdot L $ defined as follows
    	\[\label{eq:sum}
    		H_{n ,\lambda}  \cdot L  =\sum_{f \in H_{n ,\lambda}} f \cdot L.
    	\]
    	
    	The above sum \eqref{eq:sum} is finite, since we can take the sum over an $R$-basis of $H_{n, \lambda}$\footnote{$H_{n,\lambda}$ is a torsion-free module over the discrete valuation ring $R$, so it is free and has a finite $R$-basis.}. Notice also that, since $\id \in H_{n,\lambda}$ we always have $L\subset H_{n,\lambda} \cdot L$.    		
	
	 	   \subsection{Some auxiliary results from representation theory.}
	 	   Given the standard basis of $S_\lambda(V)$, the representation $\rho_{n,\lambda}$ maps a matrix $g \in \GL(n,K)$ to a matrix $\rho_{n,\lambda}(g) \in \GL(N,K)$ wih $N = \dim_K(S_\lambda(V))$; the entries of the matrix $\rho_{n,\lambda}(g)$ are homogeneous polynomials of degree $d$ in the entries of $g$. 
        
            \begin{example}\label{example:n=2_d=2}
        		Suppose that $n = d = 2$ and $\lambda = 
                            	    \ytableausetup{centertableaux}
                                    \begin{ytableau}
                                            1 & 2  \\
                                    \end{ytableau}$. Then the space $S_{\lambda}(V_{2})$ is then the second symmetric power of $V = K \cdot x_1 \ \oplus \ K \cdot x_2$ i.e. the space of homogeneous degree $2$ polynomials in $2$ variables
        		\[
        			S_{\lambda}(V) = K \cdot x_1^{\otimes 2} \ \oplus \ K \cdot (x_1 \otimes x_2 + x_2 \otimes x_1) \ \oplus \ K \cdot x_2^{\otimes 2}.
        		\]
        		The representation $\rho_{2,\lambda}$ can then be described in a matrix form as follows
        		\[
        			\rho_{2, \lambda} \colon \begin{bmatrix} a & b \\ c & d   \end{bmatrix} \mapsto \begin{bmatrix} a^2 & 2ac & c^2 \\ ab & ad + bc & cd \\ c^2 & 2bd & d^2 \end{bmatrix}. 
        		\]
        		The ring $H_{2, \lambda}$, in matrix form, given by
        		\[
        			H_{2, \lambda} \coloneqq \left\{ X \in R^{ 3 \times 3} \colon X_{12}, X_{32} \in 2 R \right\}.
        		\]
        	\end{example}

    	Let us recall the following result from representation theory of algebraic groups.
    	
    	\begin{theorem} \label{thm:PeterWeyl}
    	    Let $F$ be a field and $V$ a vector space over $F$. Suppose $\lambda$ is $\carac(F)$-core. Then the Weyl module $S_\lambda(V)$ is absolutely irreducible as a representation of the group $\GL(n, F)$. Moreover, the map
        	\begin{align} \label{eq:SchurDuality}
                    \Phi \colon S_\lambda(V) \boxtimes S_\lambda(V)^\ast &\to \Ocal(\GL(n,K))\\
        	                                            v \otimes \beta &\longmapsto \big( g \mapsto \langle \beta , \ \rho_{n,\lambda}(g) \ v \rangle   \big)   ,  \nonumber      
                \end{align}
        	where $S_\lambda(V)^\ast$ is the dual space of $S_\lambda(V)$ and $\Ocal(\GL(n,K))$ is the Hopf algebra of regular functions on $\GL(n, K)$, is injective.
    	\end{theorem}
        \begin{proof}
             When $\carac(F) = 0$ the result can be immediately deduced from the \emph{algebraic Peter-Weyl theorem},  which states that for a reductive algebraic group $G$, the ring $\Ocal(G)$ of regular functions $G$ decomposes as follows:
             \[
                \Ocal(G) = \bigoplus_{(\pi, W)\text{ irrep of } G} W \boxtimes W^\ast,
             \]
             where the sum is over the irreducible representations $\pi$ of $G$. For a reference, see \cite[Theorem 27.3.9]{TauvelRupert05} or \cite[Theorem 12.1.4]{GoodmanWallach98}.
             
             The positive characteristic case is known to the experts, but we could not find a precise reference. Hence for completeness, we provide a concise proof. Let $G \coloneqq \GL(n,F)$ and assume that  $\lambda$ is $\carac(F)$-core and denote by $W_\lambda$ the Specht module over $F$ associated to $\lambda$ (this is a representation of the symmetric group $S_d$). Since  $\lambda$ is a $\carac(F)$-core partition, the Specht module $W_\lambda$ is an absolutely irreducible and projective $S_d$-module, see \cite[Sections 11 and 23]{JamesGordon}.
             Then by the Schur-Weyl duality the Weyl module $S_\lambda(V) \cong \Hom_{S_d}(W_\lambda, V^{\otimes d})$ is an absolutely irreducible representation of $G$. So the division algebra $D := \End_{G}(S_\lambda(V))$ of intertwining  operators is trivial i.e. $D \cong F$ and using the Jacobson density theorem \cite[Theorem 3.2]{Lang}, we deduce that the $F$-linear map
             \[
                F[G] \twoheadrightarrow \End_{D}(S_\lambda(V)) = \End_{F}(S_\lambda(V)), \quad g \mapsto \rho_{n,\lambda}(g)
             \]
             is surjective, where $F[G]$ is the free group algebra of $G$ over $F$. This implies that if $\alpha \in \End_F(S_\lambda(V))^\ast$ is a linear form with $\alpha(\rho_{n,\lambda}(g)) = 0$ for all $g \in G$, then $\alpha(f) = 0$ for $f \in \End_F(S_\lambda(V)$ i.e. $\alpha = 0$. We then conclude that the map $\Phi$ is indeed injective.
        \end{proof}
    	
    	\begin{lemma}\label{lem:GL_nR_To_GL_nK}
    	    Let $A$ be a principal ideal domain with infinitely many units and $F = \Frac(A)$ its field of fractions. Let $\varphi \in \Ocal(\GL(n,F))$ be a polynomial function on $\GL(n,K)$. If $\varphi(g) = 0$ for all $g \in \GL(n,A)$, then $\varphi = 0$.
    	\end{lemma}
    	\begin{proof}
    	    Suppose that $\varphi(h) = 0$ for any $h \in \GL(n,A)$ and that $\varphi \neq 0$. Then there exists $g \in \GL(n,F)$ with $\varphi(g) = 0$. Since $A$ is a principal ideal domain, we can write the Smith normal form $g = u g' v$ of $g$ where $u,v \in \GL(n,A)$ and $g'$ is diagonal. So, by replacing $\varphi$ with the polynomial $x \mapsto \varphi(u x v)$, we may assume, without loss of generality, that $g$ is a diagonal matrix. Since $\varphi(g) \neq 0$ and $g$ is diagonal, the restriction of $\varphi$ to the space of diagonal matrices is a nonzero polynomial
    	    $$ \phi(z_1, \dots, z_n) = \varphi(\diag(z_{1}, \dots, z_{n})) = \sum_{ \substack{\nu \in \N^n} } c_\nu \ z_{1}^{\nu_1} \dots z_{n}^{\nu_n}.$$
    	    Since $\varphi$ vanishes on $\GL(n,A)$ we deduce that
    	    $$ \varphi(\diag(u_{1}, \dots, u_{n})) = \phi(u_1, \dots, u_n) = 0 ,\quad \text{for } u_1, \dots, u_n \in A^{\times}.$$
    	    Since $A$ has infinitely many units, we deduce that $\phi = 0$ which is a contradiction. So, as desired, we conclude that $\varphi = 0$.
	    \end{proof}
    	    
    	\begin{remark}
                \begin{enumerate}
                    \item The statement of Lemma \ref{lem:GL_nR_To_GL_nK} clearly fails when $A$ is a domain with finitely many units. For example if $R = \Z$, the function $g \mapsto (\det(g) - 1)(\det(g) + 1)$ vanishes on $\GL(n,\Z)$ but not on $\GL(n,\Q)$.

                    \item The ring $R$ is a discrete valuation ring. In particular, it is a principal ideal domain and has infinitely many units and the statement of Lemma \ref{lem:GL_nR_To_GL_nK} applies in this case.
                \end{enumerate}
    	\end{remark}
    	
        \begin{proposition}\label{prop:HisAnOrder}
            If $\lambda$ is $\ell$-core, the $R$-module $H_{n, \lambda}$ spans  $\End_K(S_\lambda(V))$ over $K$.
        \end{proposition}
        \begin{proof}
              Suppose that $\lambda$ is $\ell$-core and assume that $H_{n, \lambda}$ does not span $\End_K(S_\lambda(V))$. Then there exists a non-zero linear form $\alpha \in \End_K(S_\lambda(V))^\ast$ such that
              \begin{equation}\label{eq:linear_relation}
                \alpha(f) = 0, \quad \text{for all } f \in H_{n, \lambda}.  
              \end{equation}
              Equivalently, we then have $\alpha(\rho_{n,\lambda}(g)) = 0$ for all $g \in \GL(n,R)$. Since $R$ is a principal ideal domain with infinitely many units, by virtue of Lemma \ref{lem:GL_nR_To_GL_nK} we deduce that
             \[
                \alpha(\rho_{n,\lambda}(g)) = 0, \quad \text{for all } g \in \GL(n,K). 
             \]
             Here we are extending $\rho_{n,\lambda}$ in the obvious way and using the fact that $\GL(n,R)$ is an open set in $\GL(n,K)$. 
             Identifying $\End_K(S_\lambda(V))^\ast$ with $S_\lambda(V) \otimes S_\lambda(V)^\ast$ in the canonical way, \eqref{eq:linear_relation} can be rewritten as $\Phi(\alpha) = 0$. 
             But, thanks to Theorem \ref{thm:PeterWeyl}, the map $\Phi$ is injective, so we deduce that $\alpha = 0$ which is a contradiction. Hence $H_{n, \lambda}$ spans $\End(S_\lambda(V))$ over~$K$.
        \end{proof}
        
        \begin{remark}
            In the language of number theorists, when $H_{n ,\lambda}$ spans $\End_K(S_\lambda(V))$, we say that $H_{n, \lambda}$ is an \emph{order} in $\End_K(S_\lambda(V))$ i.e. a $R$-module of full rank that is also a ring.
        \end{remark}        
        
        \subsection{Buildings.}\label{sec:buildings}
                    Lattices in $S_{\lambda}(V)$ (up to scaling) are points in the \emph{Bruhat-Tits building} $\Bcal_{n,\lambda}$ for $\PGL(S_\lambda(V))$. This is an infinite simplicial complex whose $0$-simplices are the homotethy classes of lattices in $S_\lambda(V)$. For readers not familiar with Buildings, in this section we briefly recall the lattice class model of the Bruhat-Tits building associated to reductive algebraic group $\PGL$ over $K$ (giving in particular a lattice class model for $\Bcal_{n, \lambda}$). We refer to \cite{Rousseau09} for a quick introduction to the subject or \cite{AbramenkoBrown08,Brown89} for a more detailed account.

         Let $E$ be a finite dimensional vector space over $K$ and $n = \dim_K(E)$. Let $L$ be a lattice in $E$, that is an $R$-module that generates $E$ over $K$. Since $R$ is a principal ideal domain (all ideals of $R$ are generated by a power of the uniformizer $\varpi$), the lattice $L$ is a free $R$-module of rank $n$ and we can write
            \[
                L = R e_1 \oplus \dots \oplus R e_n.
            \]
            Two lattices $L$ and $L'$ in $E$ are called equivalent if $L' = u L$ for some scalar $u \in K^{\times}$ and we denote by $[L]$ the equivalence class of $L$; that is
            \[
            [L] \coloneqq \{ u L \colon u \in K^\times \}.
            \]
            \begin{definition}[Affine Buildings]
                The \emph{affine building} $\Bcal_n(K)$ associated to $\PGL(n,K)$ is a simplicial complex whose vertices are the equivalence classes $[L]$ of lattices $L$ in $E$. A set of vertices $\{ [L_0] , \dots, [L_s] \}$ is an $s$-simplex in $\Bcal_n(K)$ if and only if there exist representatives $\widetilde{L_i} \in [L_i]$ and a permutation $\sigma$ of $\{0,1, \dots, s\}$ such that 
                \[
                    \varpi\widetilde{L}_{\sigma(0)} \subsetneq \widetilde{L}_{\sigma(s)} \subsetneq \dots  \subsetneq \widetilde{L}_{\sigma(1)}  \subsetneq \widetilde{L}_{\sigma(0)}.  
                \]
                Given a basis $B = (e_1, \dots, e_n)$ of $E$, the \emph{apartment} defined by $B$ is the set of equivalence classes $[L]$ where $L$ is of the form
                \[
                    L = \varpi^{u_1} R e_1 \oplus \dots \oplus \varpi^{u_n} R e_n, \quad \text{with } u_1, \dots, u_n \in \Z.
                \]
            \end{definition}
                

            \begin{remark}
                Notice that if $\varpi L_{0} \subsetneq L_{s} \subsetneq \dots  \subset L_{1}   \subsetneq L_{0}$ for lattices $L_0, \dots, L_s$ then we have the following sequence of nested vector spaces over the residue field $k$
                \[
                    0 \subsetneq L_{s}/\varpi L_0 \subsetneq \dots  \subset L_{1}/ \varpi L_0 \subsetneq L_{0}/\varpi L_0 \cong k^n.
                \]
                Since $\dim_k(L_0 / \varpi L_0) = n$, such a sequence has length as most $n$ so the maximal simplices in $\Bcal_n(K)$ are of dimension $n-1$.
            \end{remark}

            There is a notion of min and max convexity on the Bruhat-Tits building $\Bcal_{n}(K)$. This notion was originally introduced by Faltings to study toroidal resolutions \cite{Faltings01} and it appears in the context of Mustafin varieties \cite{CHSW11,HahnLi21} (see also \cite[Section 2.1]{Zhang21}).
              \begin{figure}[H]
                \label{fig:example_d=2}
                
                        \begin{center}
                            \scalebox{0.70}{
                            \tikzset{every picture/.style={line width = 0.4pt}}          
                
                            \begin{tikzpicture}[x=0.70pt,y=0.70pt,yscale=-1,xscale=1]
                            
                             \draw  [fill={rgb, 255:red, 0; green, 0; blue, 255 }  ,fill opacity=1 ] (268.81,195.84) .. controls (268.84,193.52) and (270.62,191.67) .. (272.78,191.7) .. controls (274.94,191.73) and (276.67,193.64) .. (276.65,195.96) .. controls (276.62,198.27) and (274.84,200.13) .. (272.68,200.1) .. controls (270.51,200.07) and (268.78,198.16) .. (268.81,195.84) -- cycle ;
                             
                             \draw  [fill={rgb, 255:red, 0; green, 0; blue, 255 }  ,fill opacity=1 ] (299.03,174.06) .. controls (299.05,171.74) and (300.83,169.88) .. (303,169.91) .. controls (305.16,169.95) and (306.89,171.85) .. (306.86,174.17) .. controls (306.83,176.49) and (305.05,178.34) .. (302.89,178.31) .. controls (300.73,178.28) and (299,176.38) .. (299.03,174.06) -- cycle ;
                             
                             \draw  [fill={rgb, 255:red, 0; green, 0; blue, 255 }  ,fill opacity=1 ] (243.01,173.33) .. controls (243.04,171.01) and (244.81,169.15) .. (246.98,169.19) .. controls (249.14,169.22) and (250.87,171.12) .. (250.84,173.44) .. controls (250.81,175.76) and (249.04,177.62) .. (246.87,177.58) .. controls (244.71,177.55) and (242.98,175.65) .. (243.01,173.33) -- cycle ;
                             
                             \draw  [fill={rgb, 255:red, 0; green, 0; blue, 255 }  ,fill opacity=1 ] (210.57,183.42) .. controls (210.6,181.2) and (212.3,179.43) .. (214.36,179.46) .. controls (216.43,179.49) and (218.08,181.31) .. (218.05,183.52) .. controls (218.03,185.74) and (216.33,187.51) .. (214.26,187.48) .. controls (212.2,187.45) and (210.54,185.63) .. (210.57,183.42) -- cycle ;
                         
                             \draw   (252.02,261.95) .. controls (252.05,259.64) and (253.82,257.78) .. (255.99,257.81) .. controls (258.15,257.84) and (259.88,259.75) .. (259.85,262.07) .. controls (259.82,264.39) and (258.05,266.24) .. (255.88,266.21) .. controls (253.72,266.18) and (251.99,264.27) .. (252.02,261.95) -- cycle ;
                            
                             \draw   (286.1,261.53) .. controls (286.13,259.22) and (287.9,257.36) .. (290.07,257.39) .. controls (292.23,257.42) and (293.96,259.33) .. (293.93,261.65) .. controls (293.9,263.97) and (292.13,265.82) .. (289.96,265.79) .. controls (287.8,265.76) and (286.07,263.85) .. (286.1,261.53) -- cycle ;
                        
                            \draw [color={rgb, 255:red, 255; green, 0; blue, 0 }  ,draw opacity=1 ][line width=1.5]    (250.01,176.52) -- (269.59,193.32) ;
                        
                            \draw [color={rgb, 255:red, 255; green, 0; blue, 0 }  ,draw opacity=1 ][line width=1.5]    (276.14,193.21) -- (299.81,176.58) ;
                        
                            \draw [color={rgb, 255:red, 254; green, 19; blue, 19 }  ,draw opacity=1 ][line width=1.5]    (272.68,200.1) -- (272.9,225.53) ;
                        
                            \draw [color={rgb, 255:red, 255; green, 0; blue, 0 }  ,draw opacity=1 ][line width=1.5]    (218.05,183.52) -- (243.4,175.85) ;
                        
                            \draw    (258.01,258.29) -- (270.43,232.86) ;
                        
                            \draw [color={rgb, 255:red, 255; green, 0; blue, 0 }  ,draw opacity=1 ][line width=1.5]    (245.02,170.03) -- (227.06,146.98) ;
                        
                            \draw    (287.78,257.87) -- (275.52,232.44) ;
                        
                            \draw [color={rgb, 255:red, 255; green, 0; blue, 0 }  ,draw opacity=1 ][line width=1.5]    (303,169.91) -- (310.68,142.63) ;
                        
                            \draw   (332.99,112.71) .. controls (331.7,114.59) and (329.22,115.06) .. (327.45,113.76) .. controls (325.68,112.45) and (325.29,109.86) .. (326.58,107.98) .. controls (327.87,106.09) and (330.35,105.63) .. (332.12,106.93) .. controls (333.89,108.24) and (334.28,110.83) .. (332.99,112.71) -- cycle ;
                        
                            \draw   (364.57,96.26) .. controls (363.28,98.15) and (360.8,98.62) .. (359.03,97.31) .. controls (357.26,96) and (356.87,93.42) .. (358.16,91.53) .. controls (359.45,89.65) and (361.93,89.18) .. (363.7,90.49) .. controls (365.47,91.79) and (365.86,94.38) .. (364.57,96.26) -- cycle ;
                        
                            \draw   (336.2,76.43) .. controls (334.91,78.31) and (332.43,78.78) .. (330.66,77.47) .. controls (328.89,76.16) and (328.5,73.58) .. (329.79,71.69) .. controls (331.08,69.81) and (333.56,69.34) .. (335.33,70.65) .. controls (337.1,71.95) and (337.49,74.54) .. (336.2,76.43) -- cycle ;
                        
                            \draw    (313.76,134.72) -- (327.45,113.76) ;
                        
                            \draw    (357.62,95.72) -- (333.49,109.2) ;
                        
                            \draw    (332.8,78.43) -- (329.06,106.53) ;
                        
                            \draw   (233.55,111.15) .. controls (233.21,113.44) and (231.2,115.02) .. (229.06,114.68) .. controls (226.92,114.34) and (225.47,112.2) .. (225.81,109.91) .. controls (226.15,107.61) and (228.16,106.03) .. (230.3,106.38) .. controls (232.44,106.72) and (233.89,108.86) .. (233.55,111.15) -- cycle ;
                        
                            \draw   (254.65,81.62) .. controls (254.31,83.91) and (252.3,85.49) .. (250.16,85.15) .. controls (248.02,84.81) and (246.56,82.67) .. (246.91,80.38) .. controls (247.25,78.09) and (249.26,76.51) .. (251.39,76.85) .. controls (253.53,77.19) and (254.99,79.33) .. (254.65,81.62) -- cycle ;
                        
                            \draw   (220.84,77.12) .. controls (220.5,79.41) and (218.49,80.99) .. (216.36,80.65) .. controls (214.22,80.3) and (212.76,78.17) .. (213.1,75.87) .. controls (213.44,73.58) and (215.45,72) .. (217.59,72.34) .. controls (219.73,72.69) and (221.18,74.82) .. (220.84,77.12) -- cycle ;
                        
                            \draw    (225.85,139.89) -- (229.06,114.68) ;
                        
                            \draw    (248.22,84.39) -- (232.49,107.77) ;
                        
                            \draw    (218.68,80.5) -- (227.4,107.45) ;
                        
                            \draw   (199.13,123.72) .. controls (201,124.99) and (201.55,127.54) .. (200.36,129.41) .. controls (199.16,131.28) and (196.68,131.76) .. (194.81,130.48) .. controls (192.94,129.2) and (192.39,126.65) .. (193.58,124.79) .. controls (194.77,122.92) and (197.26,122.44) .. (199.13,123.72) -- cycle ;
                        
                            \draw   (163.94,121.48) .. controls (165.81,122.76) and (166.36,125.31) .. (165.17,127.18) .. controls (163.98,129.04) and (161.5,129.52) .. (159.63,128.25) .. controls (157.75,126.97) and (157.2,124.42) .. (158.4,122.55) .. controls (159.59,120.69) and (162.07,120.21) .. (163.94,121.48) -- cycle ;
                        
                            \draw    (221.15,142.94) -- (200.36,129.41) ;
                        
                            \draw    (181.74,98.76) -- (195.72,123.3) ;
                        
                            \draw    (166.02,124.93) -- (193.31,127.96) ;
                        
                            \draw   (185.12,94.09) .. controls (185.4,96.39) and (183.89,98.48) .. (181.74,98.76) .. controls (179.6,99.03) and (177.63,97.39) .. (177.35,95.09) .. controls (177.07,92.79) and (178.58,90.7) .. (180.73,90.42) .. controls (182.88,90.15) and (184.84,91.79) .. (185.12,94.09) -- cycle ;
                        
                            \draw  [fill={rgb, 255:red, 0; green, 0; blue, 255 }  ,fill opacity=1 ] (183.88,171.14) .. controls (186.04,171.76) and (187.32,174.01) .. (186.74,176.17) .. controls (186.16,178.32) and (183.93,179.57) .. (181.78,178.95) .. controls (179.62,178.33) and (178.34,176.08) .. (178.92,173.92) .. controls (179.5,171.76) and (181.72,170.52) .. (183.88,171.14) -- cycle ;
                        
                            \draw   (158.04,145.97) .. controls (160.2,146.59) and (161.48,148.84) .. (160.9,151) .. controls (160.32,153.16) and (158.09,154.4) .. (155.94,153.78) .. controls (153.78,153.16) and (152.5,150.91) .. (153.08,148.75) .. controls (153.66,146.6) and (155.88,145.35) .. (158.04,145.97) -- cycle ;
                        
                            \draw   (149.73,180.19) .. controls (151.89,180.81) and (153.17,183.06) .. (152.59,185.22) .. controls (152.01,187.37) and (149.79,188.62) .. (147.63,188) .. controls (145.47,187.38) and (144.19,185.13) .. (144.77,182.97) .. controls (145.35,180.81) and (147.57,179.57) .. (149.73,180.19) -- cycle ;
                            \draw [color={rgb, 255:red, 255; green, 0; blue, 0 }  ,draw opacity=1 ][line width=1.5]    (210.54,182.44) -- (186.74,176.17) ;
                            \draw    (159.94,152.91) -- (180.51,171.82) ;
                            \draw    (152.73,182.81) -- (179.61,177.03) ;
                            \draw  [fill={rgb, 255:red, 0; green, 0; blue, 255 }  ,fill opacity=1 ] (189.34,204.1) .. controls (191.01,202.55) and (193.53,202.65) .. (194.98,204.31) .. controls (196.42,205.98) and (196.24,208.59) .. (194.58,210.14) .. controls (192.91,211.69) and (190.39,211.6) .. (188.95,209.93) .. controls (187.5,208.27) and (187.68,205.66) .. (189.34,204.1) -- cycle ;
                            \draw [color={rgb, 255:red, 255; green, 0; blue, 0 }  ,draw opacity=1 ][line width=1.5]    (212.87,186.89) -- (194.98,204.31) ;
                            \draw [color={rgb, 255:red, 255; green, 0; blue, 0 }  ,draw opacity=1 ][line width=1.5]    (161.63,215.23) -- (188.1,207.41) ;
                            \draw [color={rgb, 255:red, 255; green, 0; blue, 0 }  ,draw opacity=1 ][line width=1.5]    (182.13,237.58) -- (191.85,211) ;
                            \draw   (303.33,106.18) .. controls (304.15,108.35) and (303.18,110.75) .. (301.16,111.55) .. controls (299.15,112.34) and (296.85,111.23) .. (296.03,109.06) .. controls (295.21,106.9) and (296.18,104.5) .. (298.2,103.7) .. controls (300.21,102.91) and (302.51,104.02) .. (303.33,106.18) -- cycle ;
                            \draw   (307.24,69.74) .. controls (308.06,71.91) and (307.09,74.31) .. (305.07,75.1) .. controls (303.06,75.9) and (300.76,74.79) .. (299.94,72.62) .. controls (299.12,70.45) and (300.09,68.05) .. (302.11,67.26) .. controls (304.13,66.46) and (306.42,67.57) .. (307.24,69.74) -- cycle ;
                            \draw   (276.22,82.53) .. controls (277.01,84.71) and (276.01,87.1) .. (273.99,87.87) .. controls (271.96,88.64) and (269.68,87.5) .. (268.88,85.32) .. controls (268.09,83.15) and (269.09,80.76) .. (271.12,79.99) .. controls (273.14,79.22) and (275.43,80.36) .. (276.22,82.53) -- cycle ;
                        
                            \draw    (310.63,135.09) -- (301.16,111.55) ;
                            \draw    (303.02,75.43) -- (300.78,103.79) ;
                            \draw    (275.58,87.16) -- (296.21,106.12) ;
                            \draw  [fill={rgb, 255:red, 0; green, 0; blue, 255 }  ,fill opacity=1 ] (334.81,179.42) .. controls (332.57,179.36) and (330.81,177.49) .. (330.87,175.25) .. controls (330.93,173.02) and (332.79,171.25) .. (335.03,171.32) .. controls (337.27,171.38) and (339.04,173.24) .. (338.98,175.48) .. controls (338.92,177.72) and (337.05,179.48) .. (334.81,179.42) -- cycle ;
                            \draw  [fill={rgb, 255:red, 0; green, 0; blue, 255 }  ,fill opacity=1 ] (365.21,162.11) .. controls (362.97,162.04) and (361.2,160.18) .. (361.26,157.94) .. controls (361.32,155.7) and (363.19,153.94) .. (365.43,154) .. controls (367.67,154.06) and (369.44,155.93) .. (369.38,158.17) .. controls (369.32,160.4) and (367.45,162.17) .. (365.21,162.11) -- cycle ;
                            \draw [color={rgb, 255:red, 255; green, 0; blue, 0 }  ,draw opacity=1 ][line width=1.5]    (306.3,175.12) -- (330.87,175.25) ;
                            \draw [color={rgb, 255:red, 255; green, 0; blue, 0 }  ,draw opacity=1 ][line width=1.5]    (362.3,191.13) -- (337.92,177.91) ;
                            \draw [color={rgb, 255:red, 255; green, 0; blue, 0 }  ,draw opacity=1 ][line width=1.5]    (362.32,160.31) -- (337.59,172.64) ;
                            \draw   (426.14,188.69) .. controls (423.91,188.4) and (422.33,186.36) .. (422.6,184.14) .. controls (422.88,181.92) and (424.9,180.35) .. (427.12,180.64) .. controls (429.34,180.93) and (430.93,182.97) .. (430.65,185.19) .. controls (430.38,187.41) and (428.36,188.98) .. (426.14,188.69) -- cycle ;
                            \draw   (429.75,155.99) .. controls (427.53,155.7) and (425.95,153.66) .. (426.22,151.44) .. controls (426.49,149.22) and (428.51,147.65) .. (430.74,147.94) .. controls (432.96,148.23) and (434.54,150.27) .. (434.27,152.49) .. controls (434,154.71) and (431.97,156.28) .. (429.75,155.99) -- cycle ;
                            \draw [color={rgb, 255:red, 255; green, 0; blue, 0 }  ,draw opacity=1 ][line width=1.5]    (369.07,160.51) -- (393.5,163.13) ;
                            \draw    (423.29,182.11) -- (400.27,166.49) ;
                            \draw    (426.22,151.44) -- (400.44,161.21) ;
                            \draw  [fill={rgb, 255:red, 0; green, 0; blue, 255 }  ,fill opacity=1 ] (392.13,138.45) .. controls (390.38,139.9) and (387.87,139.66) .. (386.52,137.91) .. controls (385.17,136.16) and (385.5,133.56) .. (387.25,132.11) .. controls (389,130.67) and (391.51,130.91) .. (392.86,132.66) .. controls (394.21,134.41) and (393.88,137.01) .. (392.13,138.45) -- cycle ;
                            \draw [color={rgb, 255:red, 255; green, 0; blue, 0 }  ,draw opacity=1 ][line width=1.5]    (367.67,154.21) -- (386.52,137.91) ;
                            \draw [color={rgb, 255:red, 255; green, 9; blue, 0 }  ,draw opacity=1 ][line width=1.5]    (420.43,129.04) -- (393.57,135.23) ;
                            \draw [color={rgb, 255:red, 255; green, 0; blue, 0 }  ,draw opacity=1 ][line width=1.5]    (401.23,105.48) -- (390.02,131.42) ;
                            \draw   (394.13,209.78) .. controls (392.05,208.9) and (391.04,206.51) .. (391.86,204.44) .. controls (392.68,202.37) and (395.02,201.4) .. (397.1,202.28) .. controls (399.17,203.15) and (400.19,205.54) .. (399.37,207.62) .. controls (398.55,209.69) and (396.2,210.66) .. (394.13,209.78) -- cycle ;
                            \draw   (416.96,237.89) .. controls (414.88,237.02) and (413.87,234.62) .. (414.69,232.55) .. controls (415.51,230.48) and (417.85,229.51) .. (419.93,230.39) .. controls (422,231.26) and (423.02,233.66) .. (422.2,235.73) .. controls (421.38,237.8) and (419.03,238.77) .. (416.96,237.89) -- cycle ;
                            \draw   (428.65,207.28) .. controls (426.58,206.41) and (425.56,204.02) .. (426.38,201.95) .. controls (427.2,199.87) and (429.55,198.9) .. (431.62,199.78) .. controls (433.7,200.66) and (434.71,203.05) .. (433.89,205.12) .. controls (433.07,207.19) and (430.73,208.16) .. (428.65,207.28) -- cycle ;
                            \draw    (368.92,195.34) -- (391.86,204.44) ;
                            \draw    (415.86,230.78) -- (397.55,209.5) ;
                            \draw    (426.38,201.95) -- (399.04,204.45) ;
                            \draw   (373.71,225.92) .. controls (372.83,223.78) and (373.73,221.34) .. (375.73,220.47) .. controls (377.72,219.6) and (380.04,220.62) .. (380.92,222.75) .. controls (381.8,224.89) and (380.89,227.33) .. (378.9,228.2) .. controls (376.91,229.07) and (374.58,228.05) .. (373.71,225.92) -- cycle ;
                            \draw   (370.79,262.52) .. controls (369.91,260.38) and (370.82,257.94) .. (372.81,257.07) .. controls (374.8,256.2) and (377.13,257.22) .. (378,259.35) .. controls (378.88,261.49) and (377.97,263.93) .. (375.98,264.8) .. controls (373.99,265.68) and (371.66,264.65) .. (370.79,262.52) -- cycle ;
                            \draw   (399.91,249.31) .. controls (399.04,247.18) and (399.94,244.74) .. (401.94,243.87) .. controls (403.93,242.99) and (406.25,244.01) .. (407.13,246.15) .. controls (408,248.28) and (407.1,250.72) .. (405.11,251.6) .. controls (403.12,252.47) and (400.79,251.45) .. (399.91,249.31) -- cycle ;
                            \draw    (365.63,197.29) -- (375.73,220.47) ;
                            \draw    (374.85,256.66) -- (376.32,228.21) ;
                            \draw    (401.93,243.87) -- (380.82,225.71) ;
                            \draw   (226.47,277.06) .. controls (228.36,275.82) and (230.83,276.35) .. (231.98,278.25) .. controls (233.14,280.14) and (232.53,282.68) .. (230.64,283.92) .. controls (228.74,285.16) and (226.27,284.62) .. (225.12,282.73) .. controls (223.97,280.83) and (224.57,278.29) .. (226.47,277.06) -- cycle ;
                            \draw   (191.09,279.93) .. controls (192.99,278.7) and (195.46,279.23) .. (196.61,281.13) .. controls (197.77,283.02) and (197.16,285.56) .. (195.27,286.8) .. controls (193.37,288.04) and (190.9,287.5) .. (189.75,285.61) .. controls (188.6,283.71) and (189.2,281.17) .. (191.09,279.93) -- cycle ;
                            \draw   (209.98,309.31) .. controls (211.87,308.07) and (214.34,308.6) .. (215.5,310.5) .. controls (216.65,312.39) and (216.04,314.93) .. (214.15,316.17) .. controls (212.25,317.41) and (209.78,316.87) .. (208.63,314.98) .. controls (207.48,313.08) and (208.08,310.54) .. (209.98,309.31) -- cycle ;
                            \draw    (252.47,264.22) -- (231.98,278.25) ;
                            \draw    (197.33,283.14) -- (224.69,280.1) ;
                            \draw    (213.87,308.77) -- (227.81,284.29) ;
                            \draw   (253.71,295.04) .. controls (253.72,292.72) and (255.48,290.85) .. (257.65,290.86) .. controls (259.81,290.87) and (261.56,292.76) .. (261.55,295.08) .. controls (261.54,297.4) and (259.78,299.27) .. (257.61,299.26) .. controls (255.45,299.25) and (253.7,297.36) .. (253.71,295.04) -- cycle ;
                            \draw   (237.07,327.47) .. controls (237.08,325.15) and (238.84,323.28) .. (241,323.29) .. controls (243.17,323.3) and (244.91,325.19) .. (244.9,327.51) .. controls (244.9,329.83) and (243.13,331.71) .. (240.97,331.69) .. controls (238.81,331.68) and (237.06,329.79) .. (237.07,327.47) -- cycle ;
                            \draw   (271.15,326.75) .. controls (271.16,324.43) and (272.92,322.56) .. (275.08,322.57) .. controls (277.24,322.58) and (278.99,324.47) .. (278.98,326.79) .. controls (278.97,329.11) and (277.21,330.98) .. (275.05,330.97) .. controls (272.88,330.96) and (271.14,329.07) .. (271.15,326.75) -- cycle ;
                        
                            \draw    (257.22,265.44) -- (257.65,290.86) ;
                        
                            \draw    (243.03,323.75) -- (255.24,298.22) ;
                        
                            \draw    (272.8,323.07) -- (260.33,297.75) ;
                        
                            \draw   (301.15,293.08) .. controls (300.19,291.04) and (301.09,288.7) .. (303.16,287.86) .. controls (305.22,287.01) and (307.67,287.97) .. (308.62,290.01) .. controls (309.58,292.05) and (308.68,294.39) .. (306.61,295.24) .. controls (304.55,296.09) and (302.1,295.12) .. (301.15,293.08) -- cycle ;
                        
                            \draw   (298.79,328.15) .. controls (297.83,326.11) and (298.73,323.77) .. (300.8,322.93) .. controls (302.86,322.08) and (305.31,323.04) .. (306.26,325.08) .. controls (307.22,327.12) and (306.32,329.46) .. (304.25,330.31) .. controls (302.19,331.16) and (299.74,330.19) .. (298.79,328.15) -- cycle ;
                        
                            \draw   (330.91,314.01) .. controls (329.96,311.97) and (330.86,309.62) .. (332.92,308.78) .. controls (334.98,307.93) and (337.43,308.9) .. (338.39,310.94) .. controls (339.34,312.98) and (338.44,315.32) .. (336.38,316.16) .. controls (334.31,317.01) and (331.87,316.04) .. (330.91,314.01) -- cycle ;
                        
                            \draw    (292.18,265.71) -- (303.16,287.85) ;
                        
                            \draw    (302.92,322.52) -- (303.92,295.27) ;
                        
                            \draw    (330.96,310.12) -- (308.57,292.84) ;
                        
                            \draw   (323.14,269.06) .. controls (320.82,268.7) and (319.21,266.7) .. (319.55,264.58) .. controls (319.88,262.47) and (322.02,261.03) .. (324.34,261.39) .. controls (326.65,261.74) and (328.26,263.74) .. (327.93,265.86) .. controls (327.6,267.98) and (325.45,269.41) .. (323.14,269.06) -- cycle ;
                        
                            \draw   (353.06,290.16) .. controls (350.75,289.8) and (349.14,287.8) .. (349.47,285.68) .. controls (349.8,283.57) and (351.95,282.14) .. (354.26,282.49) .. controls (356.58,282.84) and (358.19,284.84) .. (357.86,286.96) .. controls (357.52,289.08) and (355.38,290.51) .. (353.06,290.16) -- cycle ;
                        
                            \draw   (357.39,256.67) .. controls (355.07,256.32) and (353.46,254.32) .. (353.79,252.2) .. controls (354.13,250.08) and (356.27,248.65) .. (358.59,249.01) .. controls (360.9,249.36) and (362.51,251.36) .. (362.18,253.48) .. controls (361.85,255.6) and (359.7,257.03) .. (357.39,256.67) -- cycle ;
                        
                            \draw    (294.09,261.24) -- (319.55,264.58) ;
                        
                            \draw    (350.23,283.76) -- (326.54,268.03) ;
                        
                            \draw    (353.95,254.51) -- (326.83,262.98) ;
                        
                            \draw  [fill={rgb, 255:red, 0; green, 0; blue, 255 }  ,fill opacity=1 ] (312.79,134.82) .. controls (314.95,135.44) and (316.23,137.69) .. (315.65,139.85) .. controls (315.06,142) and (312.84,143.25) .. (310.68,142.63) .. controls (308.52,142.01) and (307.25,139.76) .. (307.83,137.6) .. controls (308.41,135.44) and (310.63,134.2) .. (312.79,134.82) -- cycle ;
                        
                            \draw  [fill={rgb, 255:red, 0; green, 0; blue, 255 }  ,fill opacity=1 ] (225.85,139.89) .. controls (228.01,140.51) and (229.29,142.76) .. (228.71,144.92) .. controls (228.13,147.08) and (225.91,148.32) .. (223.75,147.7) .. controls (221.59,147.08) and (220.31,144.83) .. (220.89,142.67) .. controls (221.47,140.51) and (223.7,139.27) .. (225.85,139.89) -- cycle ;
                        
                            \draw  [fill={rgb, 255:red, 0; green, 0; blue, 255 }  ,fill opacity=1 ] (398.46,160.35) .. controls (400.62,160.97) and (401.9,163.22) .. (401.32,165.38) .. controls (400.74,167.54) and (398.52,168.78) .. (396.36,168.16) .. controls (394.2,167.54) and (392.92,165.29) .. (393.5,163.13) .. controls (394.08,160.97) and (396.31,159.73) .. (398.46,160.35) -- cycle ;
                        
                            \draw  [fill={rgb, 255:red, 0; green, 0; blue, 255 }  ,fill opacity=1 ] (366.06,190.31) .. controls (368.22,190.93) and (369.5,193.19) .. (368.92,195.34) .. controls (368.34,197.5) and (366.12,198.75) .. (363.96,198.13) .. controls (361.8,197.5) and (360.52,195.25) .. (361.1,193.1) .. controls (361.68,190.94) and (363.9,189.69) .. (366.06,190.31) -- cycle ;
                        
                            \draw  [fill={rgb, 255:red, 0; green, 0; blue, 255 }  ,fill opacity=1 ] (272.9,225.53) .. controls (275.06,226.15) and (276.34,228.4) .. (275.75,230.55) .. controls (275.17,232.71) and (272.95,233.96) .. (270.79,233.34) .. controls (268.63,232.72) and (267.35,230.46) .. (267.94,228.31) .. controls (268.52,226.15) and (270.74,224.9) .. (272.9,225.53) -- cycle ;
                        
                            \draw  [fill={rgb, 255:red, 0; green, 0; blue, 255 }  ,fill opacity=1 ] (158.77,210.2) .. controls (160.93,210.82) and (162.21,213.07) .. (161.63,215.23) .. controls (161.05,217.39) and (158.83,218.63) .. (156.67,218.01) .. controls (154.51,217.39) and (153.23,215.14) .. (153.81,212.98) .. controls (154.39,210.82) and (156.61,209.58) .. (158.77,210.2) -- cycle ;
                        
                            \draw  [fill={rgb, 255:red, 0; green, 0; blue, 255 }  ,fill opacity=1 ] (403.34,97.67) .. controls (405.5,98.29) and (406.77,100.54) .. (406.19,102.7) .. controls (405.61,104.85) and (403.39,106.1) .. (401.23,105.48) .. controls (399.07,104.86) and (397.79,102.61) .. (398.38,100.45) .. controls (398.96,98.29) and (401.18,97.05) .. (403.34,97.67) -- cycle ;
                        
                            \draw  [fill={rgb, 255:red, 0; green, 0; blue, 255 }  ,fill opacity=1 ] (182.13,237.58) .. controls (184.29,238.2) and (185.57,240.45) .. (184.99,242.61) .. controls (184.41,244.76) and (182.19,246.01) .. (180.03,245.39) .. controls (177.87,244.77) and (176.59,242.52) .. (177.17,240.36) .. controls (177.75,238.2) and (179.97,236.96) .. (182.13,237.58) -- cycle ;
                            
                            \draw  [fill={rgb, 255:red, 0; green, 0; blue, 255 }  ,fill opacity=1 ] (425.39,126.26) .. controls (427.55,126.88) and (428.83,129.13) .. (428.25,131.28) .. controls (427.67,133.44) and (425.45,134.69) .. (423.29,134.07) .. controls (421.13,133.45) and (419.85,131.19) .. (420.43,129.04) .. controls (421.01,126.88) and (423.23,125.63) .. (425.39,126.26) -- cycle ;
                            
                            \end{tikzpicture} }
                    \end{center}
                        \caption{A convex set in the building $\Bcal_2(\Q_2)$ (the set of vertices colored in blue).}
                \end{figure}
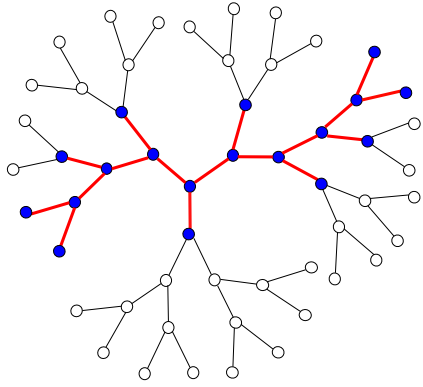   

            \begin{definition}
                A set of vertices $S = \{[L_1],\dots, [L_s]\}$ is min-convex (resp. max-convex) if for every couple of indices $1\leq i \neq j \leq s$ and representatives $\widetilde{L}_i \in [L_i]$ and $\widetilde{L}_j \in [L_j]$, we have $[\widetilde{L}_i \cap \widetilde{L}_j] \in S$ ( resp. $[\widetilde{L}_i + \widetilde{L}_j] \in S$). The set $S$ is called convex if its is both min and max convex.
            \end{definition}

            \section{Proofs of main results}\label{sec:proofs}
            \subsection{Irreducibility of Schur representations.} In this section we prove the following result, giving a sufficient condition for the irreducibility of $(\rho_{n,\lambda}, S_{\lambda}(V))$ (recall that, a representation $\rho:G\to \GL(V)$ is said to be irreducible if it has no proper invariant subspace).
             \begin{theorem}\label{thm:repIsIrred}
						If $\lambda$ is $\ell$-core, then the representation $(\rho_{n,\lambda}, S_{\lambda}(V))$ of the group $\GL(n, R)$ is irreducible.
					\end{theorem}
					 \begin{proof}
          Suppose that $\lambda$ is $\ell$-core and assume that there is a proper subspace $W \subset S_\lambda(V)$ that is $\rho_{n, \lambda}$--invariant that is
          \[
            \rho_{n, \lambda}(g)(W) \subset W, \quad \text{for all } g \in \GL(n, R).
          \]
          But, since $W$ is a proper subspace of $S_\lambda(V)$, the linear space
          \[
            \Stab(W) \coloneqq \left\{ f \in \End_K(S_\lambda(V)) \colon f(W) \subset W  \right\},
          \]
          is a proper subspace of $\End_K(S_\lambda(V))$. Since $\rho_{n,\lambda}(g) \in \Stab(W)$ we deduce that $H_{n, \lambda} \subset \Stab(W)$. By virtue of Proposition \ref{prop:HisAnOrder}, this is a contradiction. So we deduce that $\rho_{n,\lambda}$ is an irreducible representation of $\GL(n,R)$.
    \end{proof}

            \subsection{Invariant lattices: the fixed point set in the Bruhat-Tits building.}
            
             Recall that the action of $\GL(n, R)$ on $S_{\lambda}(V)$ induces an action on the building $\Bcal_{n,\lambda}$ i.e. we have a group homomorphism
			\[
				\rho_{n,\lambda}^{\Bcal} \colon \GL(n, R) \to \Aut(\Bcal_{n,\lambda}),
			\]
			sending $\GL(n, R)$ to the group of automorphisms of the building $\Bcal_{n,\lambda}$. We denote by $\Fix(\rho_{n,\lambda}^{\Bcal})$ the set of $0$-simplices of $\Bcal_{n,\lambda}$ that are fixed under $\rho_{n,\lambda}^{\Bcal}$ i.e.
			\[
				\Fix(\rho_{n,\lambda}^{\Bcal}) \coloneqq \{ [L] \colon [L] = \rho_{n,\lambda}^{\Bcal}(g)([L]) \text{ for all } g \in \GL(n, R) \}.
			\]
        Passing to the quotient modulo $\varpi$, the representation $\rho_{n,\lambda}$ of $\GL(n,R)$ naturally induces the representation
        \[
            \overline{\rho_{n,\lambda}} \colon \GL(n,k) \to \GL(S_\lambda \left( L_0 / \varpi L_0 \right)),
        \]
        where $L_0 \coloneqq R \cdot x_1 \oplus \dots \oplus R \cdot x_n$. The following result is our main theorem.

					\begin{theorem}\label{thm:FixIsOnePoint}
							Suppose that $\lambda$ is $\ell$-core. Then the set $\Fix(\rho_{n,\lambda}^{\Bcal})$ is a non-empty finite convex set in the building $\Bcal_{n,\lambda}$. Moreover, if $\lambda$ is $p$-core, the set $\Fix(\rho_{n,\lambda}^{\Bcal})$ is reduced to the one point
							\[
									\Fix(\rho_{n,\lambda}^{\Bcal}) = \left \{ [S_\lambda(L_0)] \right \},
							\]
							with $L_0 \coloneqq R \cdot x_1 \oplus \dots \oplus R \cdot x_n$.
					\end{theorem}
           
    \begin{proof}
         The lattice $L_{0,\lambda} = S_\lambda(L_0)$ is always $\rho_{n, \lambda}$-invariant so $[L_{0, \lambda}] \in \Fix(\rho_{n,\lambda}^{\Bcal})$. If $L_1, L_2$ are two $\rho_{n, \lambda}$-invariant lattices then $L_1 + L_2$ and $L_1 \cap L_2$ are also $\rho_{n, \lambda}$-invariant so we deduce that $\Fix(\rho_{n,\lambda}^{\Bcal})$ is convex in $\Bcal_{n,\lambda}$. 
         
         Now suppose that $\lambda$ is $\ell$-core. Then, by virtue of Proposition \ref{prop:HisAnOrder}, the $R$-module $H_{n, \lambda}$ is an order in $\End_K(S_\lambda(V))$. Hence there exists an integer $r > 0$ such that $\id + \varpi^r \End_R(L_{0,\lambda}) \subset H_{n ,\lambda}$. By virtue of \cite[Theorem 5.2]{ENS21}, the set of points in $\Bcal_{n, \lambda}$ that are invariant under $\id + \varpi^r \End_R(L_{0,\lambda})$ is exactly the ball $\B([L_{0, \lambda}], r)$ of center $[L_{0, \lambda}]$ and radius $r$ in the building $\Bcal_{n, \lambda}$ which is finite, so $\Fix(\rho_{n,\lambda}^{\Bcal}) \subset \B([L_{0, \lambda}], r)$ is finite.
         
         Now assume further that $\lambda$ is $p$--core and assume that there is a neighbour $[L_1]$ of $[L_{0,\lambda}]$ in $\Bcal_{n,\lambda}$; that is
         \[
            \varpi L_{0, \lambda} \subsetneq L_1 \subsetneq L_{0, \lambda},
         \]
         such that $[L_1] \in \Fix(\rho_{n,\lambda}^{\Bcal})$. Then the vector space $W = L_1 / \varpi L_{0, \lambda}$ is a proper $\overline{\rho_{n, \lambda}}$-invariant subspace of $L_{0, \lambda}/ \varpi L_{0, \lambda} = S_\lambda(L_0/\varpi L_0)$. But since $\lambda$ is $p$-core, by virtue of Theorem \ref{thm:PeterWeyl}, the representation $\overline{\rho_{n,\lambda}}$ is irreducible so this is a contradiction. We then deduce that no neighbour of $[L_{0, \lambda}]$ in $\Bcal_{n, \lambda}$ is $\rho_{n,\lambda}$-invariant. Finally, since $\Fix(\rho_{n,\lambda}^{\Bcal})$ is convex in $\Bcal_{n ,\lambda}$, we conclude that $\Fix(\rho_{n,\lambda}^{\Bcal}) = \{ [L_{0, \lambda}] \}$.
    \end{proof}
	
	\subsection{Proof of \texorpdfstring{Proposition \ref{prop:InvariantMeasures}}{}.}\label{sec:proofofcoro}

                Suppose that $\lambda$ is $\ell$-core and let $\P$ be a $\GL(n,R)$-invariant Gaussian measure on $S_\lambda(V)$ (in the sense of \cite[Section 4]{evansLocalFields2001}).
    	        Let $L \coloneqq \supp(\P)$ be the $R$-submodule of $S_\lambda(V)$ that is the support of $\P$ (see \cite[Section 4]{evansLocalFields2001}). 
    	        Then since $\P$ is $\GL(n,R)$-invariant, the lattice $L$ is $\GL(n,R)$-invariant hence also $H_{n,\lambda}$-invariant. 
    	          By virtue of Proposition \ref{prop:HisAnOrder} $H_{n,\lambda}$ is an order,  
    	        hence either $L = 0$ or $L \otimes_R K = S_\lambda(V)$. We then deduce, by virtue of Theorem \ref{thm:FixIsOnePoint}, there are only finitely $\GL(n,R)$-invariant Gaussian measures on $S_\lambda(V)$ up to scaling. 
                If moreover $\lambda$ is $p$-core, using Theorem \ref{thm:FixIsOnePoint} again, we deduce that $\P$ is the uniform distribution on a scalar multiple of the lattice $S_\lambda(L_0) \subset S_\lambda(V)$.\qed

	     \subsection{Proof of \texorpdfstring{Theorem \ref{theorem:probability}}{}}   This follows immediately from Proposition \ref{prop:InvariantMeasures} with the choice of the partition $\lambda=(d).$

	\section{Computing apartment slices of \texorpdfstring{$\Fix(\rho_{n,\lambda})$}{}}\label{sec:slices}

    In this section we assume that $\ell = 0$ and $p > 0$ (i.e. we are in the mixed characteristic case). Let $N = \dim_K(S_\lambda(V))$ and $u = (u_1, \dots, u_N)$ a basis of $S_\lambda(V)$. We denote by $\mathcal{A}_u$ the \emph{apartment} in the Bruhat-Tits building $\Bcal_{n, \lambda}$ associated to the basis $u$. We recall that the apartment $\mathcal{A}_u$ consists of the lattice classes of the form:
    \[
        \Big[ \varpi^{\alpha_1} R \ u_1 \oplus \dots \oplus \varpi^{\alpha_N} R \ u_N \Big], \quad \text{where } \alpha_1, \dots, \alpha_N \in \Z.
    \]

        We recall that an \emph{order} in $\End_K(S_\lambda(V))$ is an $R$-submodule $H$ that spans such that $H \otimes_R K = \End_K(S_\lambda(V))$ and $H$ is also a ring.
    
        \begin{definition}\label{def:graduatedOrder} We say that an order $H \subset \End_K(S_\lambda(V))$ is a \emph{graduated order} with respect to  the basis $u$ (in the sense of \cite{EHNSS22}) if there exists a matrix $M = (m_{ij}) \in \Z^{N \times N}$ with
        \begin{enumerate}
            \item $m_{ii} = 0$  for any $1 \leq i \leq N$,
            \item $m_{ij} \leq m_{ik} + m_{kj}$ for any $1 \leq i,j,k \leq N$,
        \end{enumerate}
        such that the order $H$ (in matrix form given in the basis $u$) is given by
        \[
            H = \Lambda_M \coloneqq \{ X \in K^{N \times N} \colon  X_{ij} \in \varpi^{m_{ij}} R \}.
        \]
        We denote by $\widetilde{H}_{\lambda, n}$ the smallest graduated order (with respect to $u$) containing $H_{n , \lambda}$ (see \eqref{eq:defH}) and by $\widetilde{\Fix}(\rho_{n,\lambda}^{\Bcal})$ the set of equivalence classes  of lattices in $[L]\in \Fix(\rho_{n,\lambda}^{\Bcal})$ which are fixed by the order $\widetilde{H}_{n, \lambda}$. 
        \end{definition}
        
        \begin{proposition} The following holds:
            \[
                \widetilde{\Fix}(\rho_{n,\lambda}^{\Bcal}) = \Fix(\rho_{n,\lambda}^{\Bcal})\cap \mathcal{A}_u.
            \]
        \end{proposition}
        \begin{proof}
            From \cite[Proposition 7]{EHNSS22}, we know that the fixed points of $\widetilde{H}_{n,\lambda}$ are all in the apartment $\mathcal{A}_u$. Also, since the order $\tilde{H}_{n, \lambda}$ contains the order $H_{n,\lambda}$ we deduce that the fixed points of $\tilde{H}_{n,\lambda}$ are also in $\Fix(\rho_{n,\lambda})$. 
            So we deduce that 
            \[
                \widetilde{\Fix}(\rho_{n,\lambda}^{\Bcal}) \subset (\Fix(\rho_{n,\lambda}^{\Bcal}) \cap \mathcal{A}_u).
            \]
            Now suppose that $[L] = [\varpi^{\alpha_1} R \ u_1 \oplus \dots \oplus \varpi^{\alpha_N} R \ u_N] \in \Fix(\rho_{n,\lambda}^{\Bcal}) \cap \mathcal{A}_u$. Then $[L]$ is invariant under the order $H_{n,\lambda}$ and hence
            \[
                 H_{n, \lambda} \subset \End_R(L).
            \]
            The order $\End_R(L)$ is a graduated order with respect to the basis $u$ (this is because $[L] \in \mathcal{A}_u$) so, by definition of $\tilde{H}_{n, \lambda}$, we deduce that $\widetilde{H}_{n, \lambda} \subset \End_R(L)$. So we conclude that $[L] \in \widetilde{\Fix}(H_{n, \lambda})$. This proves the other inclusion $\widetilde{\Fix}(\rho_{n,\lambda}^{\Bcal}) \supset (\Fix(\rho_{n,\lambda}^{\Bcal}) \cap \mathcal{A}_u $.
        \end{proof}

        The recipe to determine $\widetilde{\Fix}(\rho_{n,\lambda})$ from the matrix $M = (m_{i,j})$ such that $\widetilde{H}_{n,\lambda} = \Lambda_M$ is given in \cite[Proposition 7]{EHNSS22}. The set $\widetilde{\Fix}(\rho_{n,\lambda})$ is the \emph{polytrope} of points $[L_\alpha]$  where $L_\alpha \coloneqq \varpi^{\alpha_1} R \ u_1 \oplus \dots \oplus \varpi^{\alpha_N} R \ u_N $ and $\alpha$ in $\Z^{N}$ is such that:
        \[
            \alpha_i - \alpha_j \leq m_{i,j}, \quad \text{for all } 1 \leq i,j \leq N.
        \]
        Notice that the set of such $\alpha \in \Z^N$ is stable translation by the all ones vector $\bm{1} \in \Z^N$. Modulo translation by $\bm{1}$, we get a finite set $Q_M$ which is the lattice points of a \emph{polytrope} in $\R^N/\R\bm{1}$.

        \begin{example}\label{ex:computinginvariantlattices}
            Suppose that $K=\Q_2$, $n = 2$ and $\lambda = (4)$. Then the representation $\rho_{n,\lambda}$ is the action of $\GL(2,\Z_2)$ on degree $4$ homogeneous polynomials in $\Q_2[x,y]$. In the standard basis $u = (x^4, x^3y, x^2y^2, xy^3, y^4)$ this representation is given by
            \[
                \rho_{n,\lambda} \colon \begin{pmatrix}
                    a & b \\ c & d
                \end{pmatrix} \mapsto
                \begin{pmatrix}
      a^4 & a^3 c & a^2 c^2 & a c^3 & c^4\\ 
      \bm{4} a^3 b & 3a^2 b c + a^3 d & \bm{2}(a b c^2 +  a^2 c d) & b c^3+3 a c^2 d & \bm{4} c^3 d\\
      \bm{6} a^2 b^2 & \bm{3}(a b^2 c + a^2 b d) & b^2 c^2+4 a b c d+a^2 d^2 &
      \bm{3}(b c^2 d + a c d^2) & \bm{6} c^2 d^2\\
      \bm{4} a b^3 & b^3 c+3 a b^2 d & \bm{2}(b^2 c d +  a b d^2) & 3 b c d^2+a d^3 & \bm{4} c d^3\\
      b^4 & b^3 d & b^2 d^2 & b d^3 & d^4
                \end{pmatrix}
            \]
        The smallest graduated order in the basis $u$ containing $H_{n,\lambda} = \Span_R(\image(\rho_{n,\lambda}))$ is the graduated order:
        \[
            \tilde{H}_{n,\lambda} = \Lambda_{M} \coloneqq \big\{ X =(X_{ij}) \in K^{N \times N} \colon X_{ij} \in \varpi^{m_{ij}} R \big\},
        \]
        where $M$ is the exponent matrix:
        \[
            M = \begin{pmatrix}
            0 &0 &0 & 0 & 0 \\
            2 &0 &1 & 0 & 2 \\
            1 &0 &0 & 0 & 1 \\
            2 &0 &1 & 0 & 2 \\
            0 &0 &0 & 0 & 0 
            \end{pmatrix}.
        \]
        We can compute the set $Q_M$ in  $\subset \Z^5/ \Z\bm{1}$ and we get:
        \[
            Q_M = \big\{ [0,0,0,0,0], [0,1,0,1,0], [0,1,1,1,0], [0,2,1,2,0] \big\}.
        \]
        The set of invariant lattices of $\rho_{n,\lambda}$ in the standard apartment is then (up to scaling):
        
        \begin{align*}
                L_0 = & \ \Z_2 \ x^4  \oplus \ \Z_2 \ \ x^3y \oplus \Z_2 \ x^2y^2 \oplus \Z_2 \ xy^3  \oplus \Z_2 \ xy^4, \\
                L_1 =  & \  \Z_2 \ x^4 \oplus 2\Z_2 \ x^3y \oplus \Z_2 \ x^2y^2 \oplus 2\Z_2 \ xy^3 \oplus \Z_2 \ xy^4, \\
                L_2 = & \  \Z_2 \ x^4  \oplus 2\Z_2 \ x^3y \oplus 2\Z_2 \ x^2y^2 \oplus 2\Z_2 xy^3   \oplus \Z_2 \ xy^4, \\
                L_3 = & \  \Z_2 \ x^4  \oplus 4\Z_2 \ x^3y \oplus 2\Z_2 \ x^2y^2 \oplus 4\Z_2 \ xy^3 \oplus \Z_2 \ xy^4.
        \end{align*}
        \end{example}

        We now focus on the trivial partition case $\lambda = (d)$ i.e. $S_\lambda(V) = K[x_0, \dots, x_n]_{(d)}$ is the space of homogeneous polynomials of degree $d$ in $n+1$ variables. Let us denote by $u = (x^\alpha)_{|\alpha| = d}$ the standard monomial basis and $H_n = H_{n,\lambda}$ the order in $K^{N \times N}$ obtained from taking the $R$-span of the image of:
        \[
            \rho_{n,d} \colon \GL(n+1, R) \to \GL(N,R),
        \]
        with $N = \binom{n + d}{d}$. Finally denote by $\widetilde{H}_n$ the smallest graduated order with respect to the basis $u$ containing $H_n$. In this case, the computations we have conducted so far suggest the following conjecture :

        \begin{conjecture}
            The orders $\widetilde{H}_n$ and  $H_n$ are equal. Consequently, $\Fix(\rho^{\Bcal}_{n,d})$ is contained in the apartment $\mathcal{A}_u$ and can be computed as in Example \ref{ex:computinginvariantlattices}.
        \end{conjecture}

	\section{The expected number of roots of a system of invariant polynomials}\label{sec:expected}
 \subsection{\texorpdfstring{$R$}{}--structures, integrals and volumes}
 Recall that an $R$--structure on a $K$--analytic manifold $X$  is a locally constant assignment of a lattice $L_x\subset T_xX$ spanning $T_xX$ over $R$. This notion was introduced in \cite{BKL} in order to single out classical properties of nonarchimedean varieties that allow to work easily with  geometric measure theory questions.  
 
 An $R$--structure on $X$ determines a volume form $\Omega_X$, which in turn determines a Borel measure $\mu_X$ on $X$, by defining for every Borel set $U\subseteq X$ (see \cite[Section]{BKL}):
 $$\mu_X(U):=\int_{U}\Omega_X.$$
 One defined the volume of $X$ as $|X|:=\int_X\Omega_X$ (if the integral is finite).
 If $Y\subseteq X$ is an analytic submanifold, an $R$--structure on $X$ induces an $R$ structure on $Y$. 
 
Given an analytic map $\varphi:X\to Y$, where $X$ and $Y$ are $K$--analytic manifolds with $R$--strucutures, and of dimension $n$ and $m$ respectively, the \emph{absolute determinant} of $\varphi$ at $x\in X$ is defined as 
$$N(\varphi)(x):=\sigma_1\cdots \sigma_{\min{\{n,m\}}},$$
where $\sigma_1,\ldots, \sigma_{\min{\{n,m\}}}$ are the singular values of the matrix representing $D_x\varphi$ with respect to the bases given by the lattices $L_x^X$ in $T_xX$ and $L_{\varphi(x)}^Y$ in $T_{\varphi(x)}Y$ defining the $R$--structures. 
If $\varphi:X\to Y$ is injective, the volume of the image $\varphi(X)$ in $Y$ (with the induced $R$--structure) satisfies (see \cite[Theorem 3.4.4]{BKL}):
\begin{equation}\label{eq:coarea}|\varphi(X)|=\int_{X}N(\varphi)\,\Omega_X.\end{equation}
(See \cite[Section 3.4]{BKL} for more details.)
 
 In the case $X=K\textrm{P}^n$, the \emph{standard} $R$--structure is obtained as follows. For every $i=0, \ldots, n$, consider the open set 
 $$U_i:=\bigg\{[x_0, \ldots, x_n]\,\bigg|\, |x_i|=\max_j |x_j|\bigg\},$$
 together with the chart $\psi_i:U_i\to R^n$ defined by $[x_0, \ldots, x_n]\mapsto (x_0/x_i, \ldots, x_n/x_i).$ Then, for $[x]\in U_i$, we define the lattice $L_{[x]}\subset T_{[x]}K \mathrm{P}^n$ by:
 $$L_{[x]}:=D_{[x]}\psi_i^{-1}(R^n).$$
 
 In order to study the expectation of the number of solutions of random polynomial systems, we recall also the following result  from \cite{BKL} (which is a generalization of \cite[Theorem 38]{KulkarniLerario}, where this is proved in the special case $K=\Q_p$).
 
 \begin{theorem}[Projective Integral Geometry Formula]Let $\Sigma\subset \KP^m$ be a smooth analytic variety of dimension $n$. Then
 \[\label{eq:IGF}\int_{\mathrm{GL}_{m+1}(R)}\frac{|\Sigma\cap g\KP^{m-1}|}{|\KP^{n-1}|} \mathrm{d}g=\frac{|\Sigma|}{|\KP^{n}|},\]
 where ``$\mathrm{d}g$'' denotes the integration with respect to the normalized Haar measure of $\mathrm{GL}_{m+1}(R).$
 \end{theorem}
 Since the pushforward of the Haar measure of $\mathrm{GL}_{m+1}(R)$ under the map $g\mapsto g\KP^{m-1}$ induces the uniform measure on the Grassmannian $\mathbb{G}(m-1, m)$ of projective hyperplanes in $\KP^m$, \eqref{eq:IGF} can be rewritten as
 \[\label{eq:IGF2}\mathbb{E}_W\frac{|\Sigma\cap W|}{|\KP^{n-1}|} =\frac{|\Sigma|}{|\KP^{n}|},\]
 where $W\in \mathbb{G}(m-1, m)$ is sampled from the uniform distribution. Alternatively, denoting by $R^{m+1}\subset K^{m+1}$ the standard lattice, if $\zeta=(\zeta_0, \ldots,\zeta_m)$ is a point sampled from the uniform distribution on $R^{m+1}$, the random projective hyperplane $W_{\zeta}=\{\zeta_0y_0+\cdots+\zeta_my_m=0\}$ is uniformly distributed in the Grassmannian and \eqref{eq:IGF2} can be written as:
 \[\label{eq:IGF3}\mathbb{E}_{\zeta}\frac{|\Sigma\cap W_\zeta|}{|\KP^{n-1}|} =\frac{|\Sigma|}{|\KP^{n}|}.\]

 \subsection{Random system of equations}
Let now $m+1=\mathrm{dim}(K[x_0, \ldots, x_n]_{(d)})$, and consider a lattice $\Lambda \subset K[x_0, \ldots, x_n]_{(d)}$. Given a basis $\{f_j\}_{j=0, \ldots, m}$ for $\Lambda$, we define the Veronese map $\nu_{\Lambda}:K\mathrm{P}^n\to K\mathrm{P}^{m}$ associated to $\Lambda$ as: 
\begin{equation}\label{eq:nu}\nu_{\Lambda}([x_0, \ldots, x_n])=[f_0(x_0, \ldots, x_n), \ldots, f_m(x_0, \ldots, x_n)].\end{equation}
In the special case where $\Lambda$ is $\rho_{n,d}$--invariant, we have the following useful proposition.

\begin{proposition}\label{propo:constantN}Let $\Lambda\subset K[x_0, \ldots, x_n]_{(d)}$ be an invariant lattice and $\nu_\Lambda:\KP^n\to \KP^m$ be the map defined in \eqref{eq:nu}. Let $\Sigma_\Lambda:=\nu_{\Lambda}(\KP^n)$ denote the image of $\nu_{\Lambda}$ (which is a smooth submanifold) and for every point $[x]\in \KP^n$ denote by $$L_{[x]}^{\KP^n}\subset T_{[x]}\KP^{n}\quad \textrm{and}\quad L_{[y]}^{\Sigma_{\Lambda}}\subset T_{[y]}\Sigma_{\Lambda}$$  the standard $R$--structure of $\KP^{n}$ and the $R$--structure induced by the standard $R$--structure on $\KP^{\ell}$ on $\Sigma_{\Lambda}$ at the point $[y]=\nu_{\Lambda}([x])$, respectively. Then there exists $\gamma_\Lambda\in \Z$ such that for every $[x]\in \KP^{n}$
\begin{equation} \label{eq:c}D_{[x]}\nu_{\Lambda}\left(L_{[x]}^{\KP^n}\right)=p^{\gamma_{\Lambda}}\cdot L_{[y]}^{\Sigma_\Lambda}.\end{equation}
In particular, the absolute determinant of $\nu_{\Lambda}$ is a constant function:
\[\label{eq:N}N(\nu_{\Lambda})([x_0, , \ldots, x_n])\equiv N(\nu_{\Lambda})([1, 0, \ldots, 0])=p^{-n\gamma_{\Lambda}}.\]
\end{proposition}
\begin{proof}Let $g\in \mathrm{GL}_{n+1}(R)$. The fact that $\Lambda$ is $\rho_{n,d}$ invariant implies that we have a commutative diagram
$$\begin{tikzcd}
K\mathrm{P}^n \arrow[r, "\nu_{\Lambda}"] \arrow[d, "g"'] & K\mathrm{P}^{m} \arrow[d, "{\rho_{n,d}(g)}"] \\
K\mathrm{P}^n \arrow[r, "\nu_{\Lambda}"]                 & K\mathrm{P}^{m}                             
\end{tikzcd}$$
where both vertical arrows preserve the standard $R$--structures, since the image of $\rho_{n,d}$ lies in $\mathrm{GL}_{m+1}(R)$. The group $\mathrm{GL}_{n+1}(R)$ acts transitively on $K\mathrm{P}^n$ and on each tangent space, therefore it is enough to prove the first statement in the case $[x]=[1, 0, \ldots, 0]=:[e_0]$.

The stabilizer $\mathrm{Stab}_{\mathrm{GL}_{n+1}(R)}([e_0])$ can be naturally identified with $\mathrm{GL}_{n}(R)$, via the correspondence
$$\mathrm{GL}_{n}(R)\stackrel{\iota}{\hookrightarrow} \mathrm{GL}_{n+1}(R), \quad \iota:M\mapsto \left(\begin{array}{cc}1 & 0 \\0 & M\end{array}\right)$$ and the differential of $\iota(M)$ at $[e_0]$ acts on $T_{[e_0]}\KP^n\simeq \mathrm{span}\{e_1, \ldots, e_n\}\simeq K^n$ as $v\mapsto Mv.$ Let therefore $w \in L_{[e_0]}^{\KP^n}$. Then there exists $M\in \mathrm{GL}_n(R)$ such that $w=Me_1$. In order to prove \eqref{eq:c} it is enough to show that 
$$\|D_{[e_0]}\nu_{\Lambda} w\|=\|D_{[e_0]}\nu_{\Lambda} e_1\|=:p^{\gamma_{\Lambda}}.$$
To this end, using the above commutative diagram with the choice $g=\iota(M)$, we see that
\begin{align}\|D_{[e_0]}\nu_{\Lambda} w\|&=\|D_{[e_0]}\nu_{\Lambda} Me_1\|\\
&=\|D_{[e_0]}\nu_{\Lambda} (D_{[e_0]}\iota(M))e_1\|\\
&=\|D_{[e_0]}(\nu_{\Lambda} \circ \iota(M))e_1\|\\
&=\|D_{[e_0]}(\rho_{n,d}( \iota(M))\circ \nu_{\Lambda})e_1\|\\
&=\|D_{\nu_{\Lambda}([e_0])}(\rho_{n,d}( \iota(M))D_{[e_0]} \nu_{\Lambda}e_1\|\\
&=\|D_{[e_0]} \nu_{\Lambda}e_1\|,
\end{align}
where in the last step we have used the fact that $\rho_{n,d}( \iota(M))\in \mathrm{GL}_{\ell+1}(R)$, and in particular it preserves the $R$--structure. This proves \eqref{eq:c}.

 In particular, for every $[x]\in \KP^n$ the matrix representing $J(\nu_\Lambda)$ with respect to bases for the lattices giving $R$--structure in the domain and the codomain is a multiple of a matrix from $\mathrm{GL}_n(R)$ and all its singular values equal $p^{\gamma_\Lambda}$. This proves \eqref{eq:N}.
\end{proof}


\begin{definition}[The parameter of an invariant lattice]\label{def:parameter}Given a $\rho_{n,d}$--invariant lattice $\Lambda\subset K[x_0, \ldots, x_n]_{(d)}$, we define its \emph{parameter} $\delta_{\Lambda}$ by
$$\delta_{\Lambda}:=p^{-\gamma_{\Lambda}},$$
where $\gamma_\Lambda \in \Z$ is the constant given by Proposition \ref{propo:constantN}. 
\end{definition}

The next result, Theorem \ref{thm:expectation}, is the nonarchimedan analogue of \cite[Theorem 1.1]{Bu} and relates the parameter to the expected volume of random submanifolds, as in the real case.
The proof that we give here is of independent interest, compared to the proofs over the real numbers, because it shows that the integral geometry formula can be used also to deal with systems with different degrees (usually the case of different degrees is proved using the Kac--Rice formula).

Below, given a list of homogeneous polynomials $F_1, \ldots, F_c\in K[x_0, \ldots, x_n]$, we denote by $Z(F_1, \ldots, F_c)\subset \KP^{n}$ their common zero set.

\begin{theorem}\label{thm:expectation}For $i=1, \ldots, c$ let $\Lambda_i\subset K[x_0, \ldots, x_n]_{(d_i)}$ be invariant lattices and let $F=(F_1, \ldots, F_c)$ be a list of random polynomials, with the $F_i$ mutually independently sampled from the lattices $\Lambda_i$. Then
\[\label{eq:EE}\mathbb{E}_{F}\frac{|Z(F_1, \ldots F_c)|}{|\KP^{n-c}|}=\delta_{\Lambda_1}\cdots \delta_{\Lambda_c}.\]
\end{theorem}
\begin{proof}We observe first the following fact: if $Z\subset \KP^n$ is a smooth analytic variety of dimension $k$ and $\Lambda\subset K[x_0, \ldots, x_n]_{(d)}$ is an invariant lattice, then 
\[\label{eq:volZ}|\nu_{\Lambda}(Z)|=p^{-k\gamma_{\Lambda}}|Z|.\]
In fact, by \eqref{eq:c}, the map $\nu_{\Lambda}$ is a ``dilation'' of factor $p^{-\gamma_{\Lambda}}$ and the absolute determinant of $\nu_{\Lambda}|_{Z}:Z\to \nu_{\Lambda}(Z)$ is constant and equal to $p^{-k\gamma_{\Lambda}}$. The identity \eqref{eq:volZ} follows now from \eqref{eq:coarea}.

We proceed now to the proof of \eqref{eq:EE}, by induction on $c\in \N.$ If $c=1$, if $Z(F_1)\subset \KP^{n}$ is a smooth hypersurface, using \eqref{eq:volZ}, we can write
\[\label{eq:1}\frac{|Z(F_1)|}{|\KP^{n-1}|}=\frac{p^{(n-1)\gamma_{\Lambda_1}}\cdot |\nu_{\Lambda_1}(Z(F_1))|}{|\KP^{n-1}|}\]
Writing the polynomial $F_1=\sum_{|\alpha|=d_1}\zeta_\alpha x_0^{\alpha_0}\cdots x_n^{\alpha_n}$, we see that 
\[\label{eq:2}\nu_{\Lambda_1}(Z(F_1))=\nu_{\Lambda_1}(\KP^n)\cap \left\{\sum_{|\alpha|=d_1}\zeta_\alpha w_\alpha=0\right\},\]
where $\{w_\alpha\}$ are the homogeneous coordinates of $\KP^{m}=\mathrm{P}(K[x_0, \ldots, x_m]_{(d_1)})$ with respect to the monomial basis. In particular, we see that
\begin{align*}\mathbb{E}_{F_1\in \Lambda_1}\frac{|Z(F_1)|}{|\KP^{n-1}|}&=p^{(n-1)\gamma_{\Lambda_1}}\mathbb{E}_{F_1\in \Lambda_1}\frac{ |\nu_{\Lambda_1}(Z(F_1))|}{|\KP^{n-1}|}&\textrm{(using \eqref{eq:1})}\\
&=p^{(n-1)\gamma_{\Lambda_1}}\mathbb{E}_{\zeta_\alpha}\frac{\left|\nu_{\Lambda_1}(\KP^n)\cap \left\{\sum_{|\alpha|=d_1}\zeta_\alpha w_\alpha=0\right\}\right|}{|\KP^{n-1}|}&\quad \textrm{(using \eqref{eq:2})}\\
&=p^{(n-1)\gamma_{\Lambda_1}}\frac{|\nu_{\Lambda_1}(\KP^n)|}{|\KP^n|}& \textrm{(using \eqref{eq:IGF3})}\\
&=p^{(n-1)\gamma_{\Lambda_1}} \cdot p^{-n\gamma_{\Lambda_1}}&\textrm{(using \eqref{eq:N})}\\
&=p^{-\gamma_{\Lambda_1}}=\delta_{\Lambda_1}& \textrm{(by definition)}.
\end{align*}
This proves \eqref{eq:EE} in the case $c=1$. 

We proceed now with the inductive step. We assume that we have proved that
\[\label{eq:inductive}\mathbb{E}_{F}\frac{|Z(F_1, \ldots F_{c-1})|}{|\KP^{n-c+1}|}=\delta_{\Lambda_1}\cdots \delta_{\Lambda_{c-1}}\]
and we need to prove \eqref{eq:EE}.
Using \eqref{eq:1}, if $Z(F_1, \ldots, F_c)$ is smooth (which happens with probability one), we write
\begin{align}\frac{|Z(F_1, \ldots, F_c)|}{|\KP^{n-c}|}&=p^{(n-c)\gamma_{\Lambda_c}}\frac{|\nu_{\Lambda_c}(Z(F_1, \ldots, F_c))|}{|\KP^{n-c}|}\\
&=p^{(n-c)\gamma_{\Lambda_c}}\frac{|\nu_{\Lambda_c}(Z(F_1, \ldots, F_{c-1}))\cap W_{\zeta}|}{|\KP^{n-c}|},\end{align}
where $W_\zeta\simeq \KP^{m_i-1}$ is the hyperplane defined by the linear equation in the monomial variables whose coefficients are the coefficients of $F_c.$
Therefore, using \eqref{eq:IGF3}, we see that, for every $F_1, \ldots, F_{c-1}$ such that $Z(F_1, \ldots, F_{c-1})$ is smooth,
\begin{align*}\mathbb{E}_{F_c\in \Lambda_c}\frac{|Z(F_1, \ldots, F_c)|}{|\KP^{n-c}|}&=p^{(n-c)\gamma_{\Lambda_c}}\frac{|\nu_{\Lambda_c}(Z(F_1, \ldots, F_{c-1}))|}{|\KP^{n-c+1}|}\\
&=p^{(n-c)\gamma_{\Lambda_c}}\frac{p^{-(n-c+1)\gamma_{\Lambda_c}}\cdot |(Z(F_1, \ldots, F_{c-1})|}{|\KP^{n-c+1}|}\quad (\textrm{using \eqref{eq:1}})\\
&=p^{-\gamma_\Lambda}\frac{|(Z(F_1, \ldots, F_{c-1})|}{|\KP^{n-c+1}|}\\
\label{eq:*}&=\delta_{\Lambda_c}\frac{|(Z(F_1, \ldots, F_{c-1})|}{|\KP^{n-c+1}|}.\end{align*}
Taking now the expectation of \eqref{eq:*} over $F_1, \ldots, F_{c-1}$ and using the inductive assumption \eqref{eq:inductive} gives \eqref{eq:EE}. This concludes the proof.
\end{proof}
Notice in particular that from the proof of the previous theorem it follows that the parameter of and invariant lattice $\Lambda$ can also be computed as follows.
\begin{corollary}Let $\Lambda\subset K[x_0, \ldots, x_n]_{(d)}$ be an invariant lattice. Then
$$\delta_\Lambda=\mathbb{E}_{F\in \Lambda}\#(\{F=0\}\cap \KP^1).$$
\end{corollary}
\begin{proof}From the proof of the previous theorem, we have that
$$\delta_\Lambda=\mathbb{E}_{F\in \Lambda}\frac{|Z(F)|}{|\KP^{n-1}|}$$
On the other hand, integrating the previous equation, using the Integral Geometry Formula and exchanging the integrals we get
\begin{align}\delta_{\Lambda}&=\int_{\mathrm{GL}_{m+1}(R)}\mathbb{E}_{F\in \Lambda}\frac{|Z(F)\cap g\KP^{1}|}{|\KP^{0}|}\mathrm{d}g\\
&=\mathbb{E}_{F\in \Lambda}\int_{\mathrm{GL}_{m+1}(R)}\frac{|Z(F)\cap g\KP^{1}|}{|\KP^{0}|}\mathrm{d}g\\
&=\mathbb{E}_{F\in \Lambda}\#(\{F=0\}\cap \KP^1).\end{align}
\end{proof}
As a special case of Theorem \ref{thm:expectation} we also record the following corollary.
\begin{corollary}For $i=1, \ldots, n$ let $\Lambda_i\subset K[x_0, \ldots, x_n]_{(d_i)}$ be invariant lattices and let $F=(F_1, \ldots, F_n)$ be a list of random polynomials, with the $F_i$ mutually independently sampled from the lattices $\Lambda_i$. Then
\[\label{eq:EE0}\mathbb{E}_{F}\#Z(F_1, \ldots, F_n)=\delta_{\Lambda_1}\cdots \delta_{\Lambda_n}.\]
\end{corollary}

	\section{Concluding remarks and open questions} \label{sec:concluding_remarks}

    In the case where  $\lambda$ is not $\ell$-core, the set $\Fix(\rho_{n, \lambda})$ is still convex in $\Bcal_{n, \lambda}$ but can be unbounded.
	
	\begin{example}\label{exa:counterExample}
	   Assume that  $K$ is a local field of characteristic $2$ (for example the field of Laurent series $\F_2((\varpi))$). For any integer $m \geq 0$ let $L_m$ be the following lattice
	   \[
	        L_m = R \cdot x^2 \oplus \varpi^m R \cdot xy \oplus  R \cdot y^2, 
	   \]
	   in the space of homogeneous polynomial $S_{2,(2)}(V) \cong K[x,y]_{(2)}$. We claim that $L_m$ is stable under the action of $\GL(2, R)$. To see that, let $\alpha, \beta, \gamma \in R$ and  $g = \begin{pmatrix} a & b \\ c & d \end{pmatrix} \in \GL(2,R)$ and notice that:
    
	   	\begin{align*}
	   	(\alpha x^2 + \beta y^2 + \gamma \varpi^m xy) \cdot g &= \alpha (ax + by)^2 + \beta (cx + dy)^2 + \gamma \varpi^m (ax + by)(cx + dy)\\
	                                                          &= (\alpha a^2 + \beta b^2 + \gamma \varpi^m ac)x^2 + (\alpha c^2 + \beta d^2 + \gamma \varpi^m bd)y^2 \\
	                                                           &~ \ \ \ + \gamma \varpi^m (ad + bc) xy.
    	\end{align*}
	   So we deduce that $g \cdot L_m \subset L_m$ for any $g \in \GL(2,R)$. Since the action of $\GL(2,R)$ is measure preserving we deduce that $g \cdot L_m = L_m$ for any such $g$. So $\Fix(\rho_{2,(2)})$ is unbounded in $\Bcal_{2, (2)}$ since it contains the homothety class $[L_m]$ for $m \geq 0$. \footnote{Notice also that the representation $\rho_{2,(2)}$ is not irreducible (since the space W spanned by $x^2, y^2$ is $\rho_{2,(2)}$-invariant) nor semisimple (the space $W$ has no $\GL(2,K)$-invariant complement).}
	\end{example}

   	In addition to the action of $\GL(n,R)$ on the Bruhat-Tits building $\Bcal_{n, \lambda}$, one could study the action of closed subgroups of $\GL(n,R)$. The groups $\SO(n,R)$ and the symmetric group $S_n$ are of particular interest. For $\SO(n,R)$ the representation $(\rho_{n,\lambda}, S_\lambda(V))$ will no longer be irreducible and we suspect that it decomposes into a sum of irreducible representations in a similar way\footnote{Because a linear space being stable under a polynomial group action is a purely algebraic fact.} as it does over the real numbers for $\SO(n,\R)$. The set of fixed lattices for $\SO(n,R)$ will then be unbounded, as will be the case for the symmetric group $S_n$. However, it seems quite difficult to explicitly determine the set of fixed points for these groups.

\end{document}